\documentclass{gtart}

\usepackage{times}
\usepackage{amsmath, amssymb, latexsym, amscd, graphicx, epstopdf}
\usepackage{euscript}
\usepackage{subfigure}
\usepackage{epsfig}
\usepackage{psfrag}

\usepackage[inner=30mm, outer=30mm, textheight=225mm]{geometry}

\def\bkC{{\rm \kern.24em \vrule width.05em height1.4ex depth-.05ex 
\kern-.26em C}}
\def\C{\bkC}
\def\bksC{{\rm \kern.24em \vrule width.05em height1ex depth-.05ex 
\kern-.26em C}}

\def\bkE{{\rm I\kern-.22em E}}
 
\def\bkH{{\rm I\kern-.22em H}}
\def\H{\bkH} 
\def\bkN{{\rm I\kern-.17em N}}
\def\NN{\bkN}
\def\bkQ{{\rm \kern.24em \vrule width.05em height1.4ex depth-.05ex 
\kern-.26em Q}}

\def\bkR{{\rm I\kern-.17em R}}
\def\RR{\bkR}
\def\bkZ{{\rm Z\kern-.32em Z}}
\def\Z{\bkZ}
\def\bksZ{{\rm Z\kern-.22em Z}}

\def\SL{SL_2(\C)}
\def\PSL{PSL_2(\C)}

\def\D{\mathfrak{D}}
\def\Ei{\mathfrak{E}}
\def\PEi{\overline{\Ei}}

\def\PX{\overline{\mathfrak{X}}}

\def\tree{\mathfrak{T}}

\def\tri{\mathcal{T}}
\def\fol{\mathcal{F}}

\def\l{\mathcal{L}}
\def\m{\mathcal{M}}
\def\N{PF}

\def\prho{\overline{\rho}}

\DeclareMathOperator{\tr}{tr}

\DeclareMathOperator{\im}{im}

\DeclareMathOperator{\pee}{\overline{e}}

\theoremstyle{plain}

\newtheorem{thm}{Theorem}[section]
\newtheorem*{thm*}{Theorem}
\newtheorem{lem}[thm]{Lemma}
\newtheorem*{lem*}{Lemma}
\newtheorem{cor}[thm]{Corollary}
\newtheorem*{cor*}{Corollary}

\newtheorem*{cla*}{Claim}
\newtheorem{pro}[thm]{Proposition}
\newtheorem*{pro*}{Proposition}
\newtheorem{rem}[thm]{Remark}
\newtheorem*{rem*}{Remark}
\newtheorem{defn}[thm]{Definition}
\newtheorem*{defn*}{Definition}

\addtocounter{MaxMatrixCols}{4}
\begin{document}
\title{Degenerations of ideal hyperbolic triangulations} 
\author{Stephan Tillmann}

\begin{abstract}
Let $M$ be a cusped 3--manifold, and let $\tri$ be an ideal triangulation of $M.$ The deformation variety $\D(\tri),$ a subset of which parameterises (incomplete) hyperbolic structures obtained on $M$ using $\tri,$ is defined and compactified by adding certain projective classes of transversely measured singular codimension--one foliations of $M.$ This leads to a combinatorial and geometric variant of well--known constructions by Culler, Morgan and Shalen concerning the character variety of a 3--manifold.
\end{abstract}


\maketitle


\section{Introduction}

Let $M$ be the interior of an orientable, compact, connected 3--manifold with non-empty boundary consisting of a pairwise disjoint union of tori, and let $\tri$ be an ideal (topological) triangulation of $M.$ Following an approach sketched by Thurston \cite{t81}, projective classes of transversely measured singular codimension--one foliations of $M$ are associated to degenerations of ideal hyperbolic triangulations of $M,$ which are parameterised by an affine algebraic set $\D(\tri),$ called the \emph{deformation variety}. The set $\D (\tri)$ is also related to the study of representations of $\pi_1(M)$ into $\PSL,$ and used to obtain an explicit understanding of the link between the topology of a cusped 3--manifold and ideal points of varieties related to representations into $\SL$ and $\PSL.$ This link is known to exist from constructions by Culler, Morgan and Shalen~\cite{cs, ms1}.

The new contribution of this paper lies in an analysis of the deformation variety using hyperbolic geometry and tropical geometry, as well as ideas of Bestvina~\cite{best}, Paulin~\cite{P}, Thurston~\cite{t81} and Yoshida~\cite{y}. The canonical Morgan-Shalen compactification of the character variety is infinite dimensional and relies on Hironaka's theorem. The approach taken here uses an ideal triangulation as a particular, finite-dimensional coordinate system. The aforementioned link is explicitly described as a degeneration of the ideal hyperbolic triangulation, and the use of ideal triangulations and tropical geometry allows an algorithmic construction of surfaces (or, more generally, transversely measured singular foliations) dual to ideal points of the character variety.

\subsection{Ideal points via tropical geometry}

\begin{figure}[t]
\begin{center}
      \subfigure[]{\label{fig:deg singular foliation}
        \includegraphics[height=4cm]{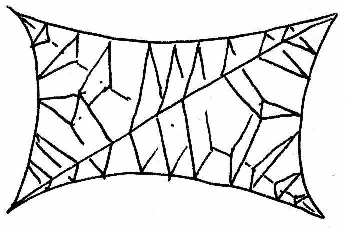}}
      \qquad\qquad
      \subfigure[]{ \label{fig:non-deg singular foliation}
        \includegraphics[height=4cm]{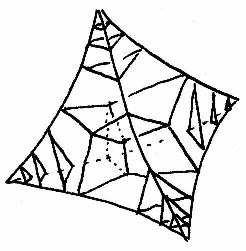}}
   \end{center}
    \caption[]{The two types of singular foliations of ideal tetrahedra}
    \label{fig:singular foliation}
\end{figure}

The set of ideal points $\D_\infty (\tri)$ is defined to be Bergman's \emph{logarithmic limit set} of $\D (\tri).$ This tropical compactification turns out to have an explicit geometric interpretation:

An ideal hyperbolic tetrahedron in $\H^3$ is the convex hull of four distinct points on the sphere at infinity. It can be positively or negatively oriented and is flat (an ideal quadrilateral) if its ideal vertices lie on a round circle. It degenerates as one of its ideal vertices is moved to coincide with another. The deformation variety describes the shapes of ideal hyperbolic tetrahedra, and has the property that an ideal point is approached if and only if some tetrahedron degenerates. 
A degenerating tetrahedron is seen to become very long and thin. It converges in the Gromov--Hausdorff sense to a dual spine if the hyperbolic metric is suitably rescaled, and it inherits a singular codimension--one foliation from the collapse as shown in Figure \ref{fig:deg singular foliation}. Applying the same rescaling process to a tetrahedron which does not degenerate gives the singular foliation shown in Figure \ref{fig:non-deg singular foliation}. Moreover, the foliation inherits a transverse measure from the relative growth rates at which the tetrahedra degenerate.

The main result in Section~\ref{sec:Ideal points and normal surfaces} (Proposition~\ref{comb:homeo}) formalises this interpretation by giving a canonical identification of $\D_\infty (\tri)$ with a compact subset of the projective admissible solution space of spun-normal surface theory, denoted $\N(\tri)$ in \cite{part1}, via a natural injection;
$$\D_\infty (\tri) \ni \xi \to N(\xi) \in \N(\tri).$$
The set $\D_\infty (\tri)$ can be computed using work of Bogart, Jensen, Speyer, Sturmfels and Thomas~\cite{bjsst}; the algorithm is implemented in the software package \tt{gfan}\rm\ by Jensen \cite{gfan}.

\subsection{A sufficient criterion for non-trivial actions}

The ideal points of $\D(\tri)$ are analysed via actions on $\RR$--trees. The first part, found in Section~\ref{sec:The dual tree}, applies to any admissible projective class. Every $N \in \N (\tri)$ defines a transversely measured singular codimension--one foliation $\fol.$ If $N$ has rational coordinate ratios, then there is also an associated spun-normal surface $S$ corresponding to a closed (not necessarily connected) leaf. The foliation lifts to a foliation $\widetilde{\fol}$ of the universal cover $\widetilde{M}$ of $M,$ and the leaf space $\widetilde{M}/\widetilde{\fol}$ is turned into an $\RR$--tree $\tree_N$ on which $\pi_1(M)$ acts by isometries with translation length function $l_N.$ 

How can one decide whether $l_N$ is non-trivial? This question is addressed in Section~\ref{The geometric construction} for the admissible projective classes corresponding to ideal points of $\D (\tri).$ So suppose that $N = N(\xi)$ for some $\xi \in \D_\infty(\tri).$ Each ideal 3--simplex in $\widetilde{M}$ is imbued with a hyperbolic structure, and $\widetilde{M}/\widetilde{\fol}$ is suitably interpreted as a limit of $\widetilde{M}$ as the ideal triangulation degenerates. There is an associated degeneration of the hyperbolic space $\H^3$ with a limiting action of $\pi_1(M)$ with length function $l_\xi$ coming from a sequence of representations. The function $l_\xi$ is not canonical, but it has the property that 
$$0 \le  l_\xi(\gamma) \le l_N(\gamma)$$ 
for all $\gamma \in \pi_1(M).$ It is in general difficult to determine $l_\xi(\gamma).$ An exception are peripheral elements, where it is determined by the linear functional $\nu_N$ of \cite{part1}, \S 3. The relationship of $\nu_N$ with the action on the tree is established in Section \ref{Holonomy variety and essential surfaces}, where it is also shown that the boundary curves of the spun-normal surfaces in $\D_\infty(\tri)$ are the boundary curves of essential surfaces which are strongly detected by the character variety.

Before the next result is stated, more terminology is introduced. A surface in $M$ is \emph{non-trivial} if it is essential or can be reduced to an essential surface by performing compressions and then possibly discarding some components. If $\xi$ has rational coordinate ratios, so does $N(\xi)$ and there is an associated spun-normal surface $S(\xi),$ which determines the foliation. A compact spun-normal surface is an ordinary normal surface. 

\begin{thm}\label{thm:main 1}
Let $M$ be the interior of a compact, connected, orientable, irreducible 3--manifold with non-empty boundary consisting of a disjoint union of tori, and let $\tri$ be an ideal triangulation of $M.$ Let $\xi \in \D_\infty(\tri)$ and $N = N(\xi) \in \N (\tri).$ 

The action on $\tree_N$ is non-trivial if $\nu_N(\gamma)\neq 0$ for some peripheral element $\gamma \in \pi_1(M).$ In particular, if $\xi$ has rational coordinate ratios and the associated surface $S(\xi)$ is non-compact, then $S(\xi)$ is non-trivial.
\end{thm}

Three further remarks to the above theorem should be added. First, it can be viewed as a generalisation of Yoshida's main result in \cite{y}. Second, it opens the door to an algorithmic approach to the character variety techniques of Culler, Morgan and Shalen in general, and to computing boundary curves strongly detected by the character variety in particular. Third, it is proved without the use of Culler-Morgan-Shalen theory, but merely with combinatorial and geometric arguments. 

Not all ideal points of $\D(\tri)$ give rise to non-trivial actions since there is a trivial, closed normal surface associated to an ideal point of the deformation variety of the Whitehead link complement (with its standard triangulation); see \cite{part4}. It remains an open problem to determine a necessary and sufficient condition purely in terms of normal surface theory for the ideal points, where $\nu_N(\gamma)= 0$ for all peripheral elements $\gamma \in \pi_1(M).$ In the next main result, a sufficient condition using Culler-Morgan-Shalen theory is given.

\subsection{Relationship with the Morgan--Shalen compactification}

Morgan and Shalen \cite{ms1} compactified the character variety of a 3--manifold by identifying ideal points with certain actions of $\pi_1(M)$ on $\RR$--trees (given as points in an infinite dimensional space), and dual to these are codimension--one measured laminations in $M.$ This paper has taken a different but related approach by compactifying the deformation variety with certain transversely measured singular codimension--one foliations (a finite union of convex rational polytopes), and dual to these are (possibly trivial) actions on $\RR$--trees. The relationship between these compactifications is investigated in Section \ref{The geometric construction}, with focus on the original consideration of ideal points of curves and surfaces dual to Bass--Serre trees due to Culler and Shalen \cite{cs}. A surface which can be reduced to an essential surface that is detected by the character variety in the sense of \cite{tillus_mut} will be called \emph{weakly dual to an ideal point of a curve in the character variety}. 

\begin{thm}\label{thm: morgan-shalen}
Let $M$ be the interior of a compact, connected, orientable 3--manifold with non-empty boundary consisting of a disjoint union of tori, and let $\tri$ be an ideal triangulation of $M.$ Let $\{Z_i\} \subset \D(\tri)$ be a sequence approaching the ideal point $\xi\in\D_\infty(\tri),$ and denote $\{\chi_i\} \subset \PX(M)$ the associated sequence of characters. Then: 
\begin{enumerate}
\item The action on $\tree_N$ is non-trivial if an ideal point of the character variety is approached by the sequence $\{\chi_i\}$ in the sense of Morgan and Shalen.
\item If $\xi$ has rational cooordinate ratios and for some $\gamma\in \pi_1(M),$ $\{ | \chi_i(\gamma)| \}$ is unbounded, then $S(\xi)$ is non-trivial and (weakly) dual to an ideal point of a curve in the character variety of $M$.
\end{enumerate}
\end{thm}

\begin{figure}[t]
\begin{center}
        \includegraphics[width=8cm]{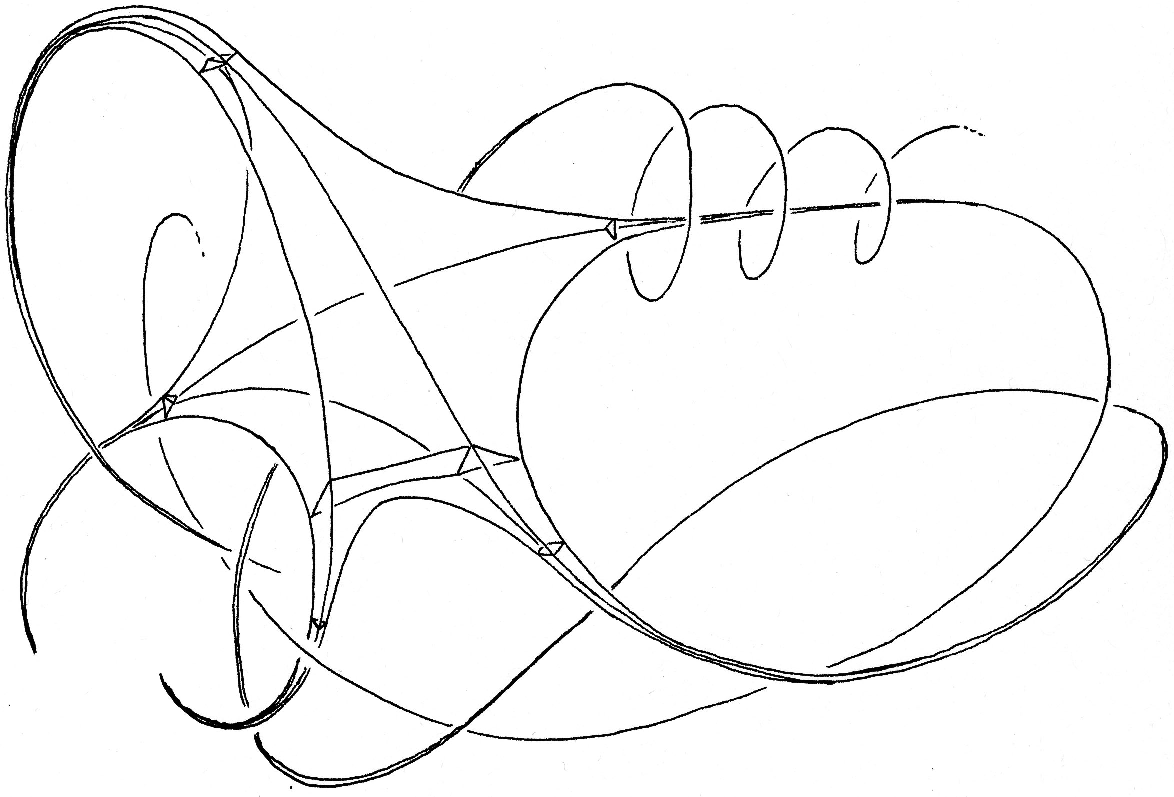}
   \end{center}
    \caption[]{A spun-normal surface develops in a degenerating ideal triangulation: shown are one tetrahedron which degenerates (to the left) and one which does not (to the right). This figure is based on Figure 4 in \cite{t1}.}
    \label{fig:degenerating triangulation}
\end{figure}

\subsection{Examples and further extensions}

The reader may refer to the example in Section~\ref{sec: figure eight} whilst reading this paper. This example also illustrates the nicest setting: an ideal point is approached through positively oriented ideal hyperbolic triangulations. The tetrahedra associated to the degeneration become very long and thin, and a surface develops in the thin part, along which the manifold will split apart. This is illustrated in Figure~\ref{fig:degenerating triangulation}. The surface naturally inherits a cell decomposition as a spun-normal surface with respect to the ideal triangulation. In this case, one obtains a splitting of the limiting Dehn--Thurston surgered 3--manifold along Euclidean (possibly cone) 2--manifolds into 3--dimensional pieces, each of which is either a hyperbolic (possibly cone) 3--manifold or a Seifert fibered manifold; the details are worked out in \cite{part3} using angle structures and standard spines.

The definition of the deformation variety applies to any triangulation of a closed 3--manifold \cite{L, LTY} or ideal triangulation of a topologically finite non--compact 3--manifold \cite{ST}, and many of the results proved in this paper go through verbatim. Examples show that the ``fake" ideal points of the deformation variety can bear interesting information; e.g.\thinspace if one applies the construction to the minimal triangulation of quaternionic space $S^3/Q_8,$ one obtains a curve with three ideal points corresponding to the three 1--sided Heegaard splitting surfaces even though the map to the character variety is constant.


\subsection*{Acknowledgements} 

The author thanks Steven Boyer, Craig Hodgson, Feng Luo and Walter Neumann for helpful discussions. This research was supported by a CRM/ISM Postdoctoral Fellowship and under the Australian Research Council's Discovery funding
scheme (project number DP110101104).


\section{The Deformation variety}
\label{The deformation variety}

In this section, notation is fixed and several facts and standard results are stated.
Throughout this paper, $M$ denotes the interior of a compact, connected, orientable 3--manifold with non-empty boundary consisting of a pairwise disjoint union of tori. Let $\tri$ be any topological ideal triangulation of $M.$ The \emph{deformation variety} $\D (\tri)$ describes hyperbolic structures for the ideal 3--simplices in $M$ subject to certain gluing equations. It is birationally equivalent to Thurston's \emph{parameter space}, which appears for example in \cite{t, nz, ch, ac}. The latter requires the choice of an edge for each tetrahedron, whilst the former keeps the symmetries of the triangulation.


\subsection{Ideal triangulations}

The ideal triangulation $\tri$ of $M$ consists of a pairwise disjoint union of standard Euclidean 3--simplices, $\widetilde{\Delta} = \cup_{k=1}^{n} \widetilde{\Delta}_k,$ together with a collection $\Phi$ of Euclidean isometries between the 2--simplices in $\widetilde{\Delta};$ termed \emph{face pairings}. Then $M = (\widetilde{\Delta} \setminus \widetilde{\Delta}^{(0)} )/ \Phi,$ and $P = \widetilde{\Delta} / \Phi$ is the associated \emph{pseudo-manifold} (or \emph{end-compactification} of $M$) with quotient map $p\co \widetilde{\Delta} \to P.$ Let $\sigma$ be a $k$--simplex in $\widetilde{\Delta}.$ Then $p(\sigma)$ may be a singular $k$--simplex in $P,$ and is termed a \emph{$k$--singlex} for short. Denote $\Sigma^k$ the set of all $k$--singlices in $P.$ An \emph{ideal $k$--simplex} is a $k$--simplex with its vertices removed. The vertices of the $k$--simplex are termed the \emph{ideal vertices} of the ideal $k$--simplex. Similar for singlices. The standard terminology of (ideal) edges, (ideal) faces and (ideal) tetrahedra will be used for the singlices in $M$ and $P.$

By hypothesis on $M,$ the link of each vertex in $P$ is a torus, so $\chi(M)=0$ and $|\Sigma^3| = |\Sigma^1|.$ These facts are irrelevant for much of the material below. However, the case of spherical vertex links is treated in \cite{L, LTY}, and for the case of higher genus boundary components, it would be more natural to develop a theory allowing hyperideal vertices in order to capture representations whose restriction to the boundary is non-abelian.


\subsection{Definitions and notation}
\label{combinatorics:Deformation variety}

Let $\Delta^3$ be the standard 3--simplex with a chosen orientation. Suppose the edges from one vertex of $\Delta^3$ are labeled by $z,$ $z'$ and $z''$ so that the opposite edges have the same labeling. Then the cyclic order of $z,$ $z'$ and $z''$ viewed from each vertex depends only on the orientation of the 3--simplex. It follows that, up to orientation preserving symmetries, there are two possible labelings, and we fix one of these labelings. 

Suppose $\Sigma^3 = \{ \sigma_1, \ldots, \sigma_n\}.$ Since $M$ is orientable, the 3--simplices in $\Sigma^3$ may be oriented coherently. For each $\sigma_i \in \Sigma^3,$ fix an orientation preserving simplicial map $f_i\co \Delta^3 \to \sigma_i.$ Let $\Sigma^1 = \{ e_1, \ldots, e_n\},$ and let $a^{(k)}_{ij}$ be the number of edges in $f_i^{-1}(e_j)$ which have label $z^{(k)}.$

For each $i \in \{1,\ldots, n\},$ define
\begin{equation}\label{eq:para}
    p_i   = z_i (1 - z''_i) - 1,\quad
    p'_i  = z'_i (1 - z_i) - 1,\quad
    p''_i = z''_i (1 - z'_i) - 1,
\end{equation}
and for each $j \in \{1,\ldots, n\},$ let
\begin{equation}\label{eq:glue}
    g_j = \prod_{i=1}^n z_i^{a_{ij}} {(z'_i)}^{a'_{ij}} {(z''_i)}^{a''_{ij}}
    -1. 
\end{equation}
Setting $p_i=p'_i=p''_i = 0$ gives the \emph{parameter equations}, and setting $g_j=0$ gives the \emph{hyperbolic gluing equations}. For a discussion and geometric interpretation of these equations, see \cite{t, nz}. The parameter equations imply that $z_i^{(k)} \neq 0,1.$

\begin{defn}
The \emph{deformation variety $\D (\tri)$} is the variety in $(\C-\{0 \})^{3n}$ defined by the parameter equations and the hyperbolic gluing equations.
\end{defn}


\subsection{Ideal hyperbolic tetrahedra}

Four ordered points on $\partial \H^3 = \C \cup \{\infty\}$ determine a cross ratio by the following formula:
\begin{equation*}
(v_i,v_j;v_k,v_l) = \frac{v_i-v_k}{v_i-v_l}\cdot\frac{v_j-v_l}{v_j-v_k}.
\end{equation*}
If the triple $(z, z', z'')$ of complex numbers satisfies the three equations
\begin{equation*}
     z (1 - z'') = 1,\quad
     z' (1 - z) = 1,\quad
     z'' (1 - z') = 1,
\end{equation*}
then there is an oriented ideal hyperbolic tetrahedron $[v_0, v_1, v_2, v_3]$ in $\H^3,$ where $v_i \in \partial \H^3$ and the order determines the orientation, such that 
\begin{equation*}
     z = (v_0, v_1; v_2, v_3),\quad
     z' = (v_0, v_3; v_1, v_2)\quad
     z'' = (v_0, v_2; v_3, v_1).
\end{equation*}
The ideal hyperbolic tetrahedron is unique up to orientation preserving isometries, and hence congruent with $[0,1,\infty, z].$ It is termed \emph{positively oriented} if $\Im(z)>0,$ \emph{flat} if $\Im(z)=0,$ and \emph{negatively oriented} if $\Im(z)<0.$


\subsection{Essential edges}

Denote $C$ a compact core of $M;$ this is obtained by removing a small open regular neighbourhood from each vertex in $P,$ such that $C$ inherits a decomposition into truncated tetrahedra. Let $e$ be an ideal edge in $M.$ The intersection $\alpha = e\cap C$ is \emph{homotopic into $\partial C$} if and only if there is (1) an arc $\beta \subset \partial C$ such that $\partial \beta = \partial \alpha$ and (2) a map $f \co D^2 \to C$ such that $f(\partial D^2) =\alpha \cup \beta$ and the restriction of $f$ to the boundary is a local homeomorphism. The edge $e$ in $M$ is \emph{essential} if $e\cap C$ is not homotopic into $\partial C,$ and it is \emph{not essential} otherwise.

For instance, if $M$ is the complement of a knot or link in $S^3,$ then $P$ is simply connected. Whence each edge is null-homotopic in $P$ but all edges may be essential in $M;$ see Figure~\ref{fig:fig8_tets} for an example. A proof of the following result can be found in \cite{Dun, ST}:

\begin{lem}
If $\D(\tri)$ is non-empty, then all edges in $M$ are essential.
\end{lem}


\subsection{End-compactification of universal cover}

Denote $p \co \widetilde{M} \to M$ the universal cover of $M,$ and lift the ideal triangulation of $M$ to a $\pi_1(M)$--equivariant ideal triangulation of $\widetilde{M}.$ Denote $C \subset M$ the compact core as above. Let $\widehat{P}$ be the space obtained from the universal cover $\widetilde{C}$ of $C$ by attaching the cone over each connected boundary component $B$ of $\widetilde{C}$ to a point $v_B.$ We then have natural inclusions 
$$\widetilde{C} \subset \widetilde{M} \subset \widehat{P},$$
and $\widehat{P}$ is termed the \emph{end-compactification of $\widetilde{M}$ with respect to $M.$} Note that $\widehat{P}$ is also simply connected since adding cones over connected spaces does not increase fundamental group, and that there is a natural, simplicial map $\widehat{P} \to P,$ where the simplicial structure arises from lifting the decomposition of $C$ into truncated tetrahedra and completing them to tetrahedra when adding the cones. The map $\widehat{P} \to P$ is also denoted by $p$ since it restricts to $p \co \widetilde{M} \to M.$

It is hoped that the notation and terminology does not lead to confusion. For instance, when $M$ is hyperbolic, then $\widetilde{M}$ is an open 3--ball and the natural compactification of this open 3--ball (without reference to $M$) is homeomorphic to the closed 3--ball, whilst the end-compactification $\widehat{P}$ is $\widetilde{M}$ with countably many points added. Also, as noted above, $P$ itself may be simply connected.


\subsection{Developing maps and characters}
\label{sec: developing}

The following facts can be found in \cite{y}, \S 5; see also \cite{ST} for a more detailed discussion.

\begin{lem}
Let $M$ be the interior of a compact, connected, orientable 3--manifold with non-empty boundary consisting of a pairwise disjoint union of tori, and let $\tri$ be an ideal triangulation of $M.$ For each $Z \in \D(\tri),$ there exists a representation $\prho_Z \co \pi_1(M) \to \PSL$ and a $\prho_Z$--equivariant, continuous map $D_Z \co \widehat{P} \to \overline{\H}^3$ with the property that for each ideal tetrahedron $\sigma \subset \widetilde{M},$ the image $D_Z(\sigma)$ is a hyperbolic ideal tetrahedron with edge invariants determined by $Z.$ Moreover, $\prho_Z$ is well-defined up to conjugacy and the well-defined map 
$$\chi_\tri \co \D(\tri) \to \PX (M)$$
is algebraic.
\end{lem}


\subsection{Holonomies and eigenvalue variety}
\label{sec: Holonomies and eigenvalue variety}

As above, let $C$ be a compact core of $M.$ Each boundary torus $T_i,$ $i=1,{\ldots} ,h,$ of $C$ inherits a triangulation $\tri_i$ induced by $\tri.$ Let $\gamma$ be a closed simplicial path on $T_i.$ In
\cite{nz}, the \emph{holonomy} $\mu(\gamma)$ is defined as $(-1)^{|\gamma|}$ times the product of the moduli of the triangle vertices touching $\gamma$ on the right, where $|\gamma|$ is the number of 1--simplices of $\gamma,$ and the moduli asise from the corresponding edge labels.

At $Z \in \D (\tri),$ evaluating $\mu(\gamma)$ gives a complex number $\mu_Z(\gamma) \in \C \setminus\{0\}.$ It is stated in \cite{t, nz}, that $\mu_Z(\gamma)$ is the square of an eigenvalue; one has: 
$$(\tr \prho_Z (\gamma))^2 = \mu_Z(\gamma) + 2 + \mu_Z(\gamma)^{-1}.$$ 
This can be seen by putting a common fixed point of $\prho_Z (\pi_1(T_i))$ at infinity in the upper--half space model, and writing
$\prho_Z(\gamma)$ as a product of M\"obius transformations, each of
which fixes an edge with one endpoint at infinity and takes one face
of a tetrahedron to another. 

Choose a basis $\{\m_i, \l_i\}$ of $\pi_1(T_i)\cong H_1(T_i)$ for each boundary torus $T_i.$  Since $\mu_Z\co \pi_1(T_i) \to \C\setminus \{0\}$ is a homomorphism for each $i=1,{\ldots} ,h,$ there is a well--defined rational map:
\begin{equation*}
    \pee \co \D (\tri) \to (\C-\{0\})^{2h} \qquad
    \pee(Z) = (\mu_Z (\m_1),{\ldots} ,\mu_Z (\l_h )).
\end{equation*}
The closure of its image is contained in the \emph{$\PSL$--eigenvalue variety} $\PEi (M)$ of \cite{tillus_ei}.


\subsection{Remarks on hyperbolic manifolds}

For cusped hyperbolic 3--manifolds, there are a number of special facts that are well-known to follow from work in \cite{t, nz}. Two are stated below to highlight the fact that the main results of this paper have interesting applications in the study of the Dehn surgery components in the character variety, i.e.\thinspace the components containing the characters of discrete and faithful representations.

\begin{thm}
Let $M$ be an orientable, connected, cusped hyperbolic 3--manifold. Let $\tri$ be an ideal triangulation of $M$ with the property that all edges are essential. Then there exists $Z \in \D(\tri),$ such that $\prho_Z\co \pi_1(M)\to \PSL$ is a discrete and faithful representation. Moreover, the whole {Dehn surgery component} $\PX_0(M)$ containing the character of $\prho_Z$ is in the image of $\chi_\tri.$
\end{thm}

As in \cite{LTY}, one can identify the discrete and faithful representation algorithmically. One first imposes the \emph{completeness equation} 
\begin{equation*}
        (\mu_Z (\m_1),{\ldots} ,\mu_Z (\l_h )) = (1,\ldots, 1),
\end{equation*}
written $\pee(Z)=1$ in short-hand, in order to ensure that all peripheral subgroups have parabolic representations. This gives a finite collection of points on $\D(\tri).$ The \emph{volume} of $Z \in \D (\tri)$ can be defined using the Lobachevsky--Milnor formula, and it follows from Francaviglia \cite{Fr} that it suffices to identify a point of maximum volume.

There is a special case, where one can start with very little topological information; this is Thurston's method to construct hyperbolic structures using ideal triangulations. Denote $\C_>$ (resp.\thinspace $\C_\ge$) the set of all complex numbers with positive (resp.\thinspace non-negative) imaginary part.

\begin{thm}
Let $M$ be the interior of a compact, connected, orientable 3--manifold with non-empty boundary consisting of a pairwise disjoint union of tori. Let $\tri$ be an arbitrary ideal triangulation of $M.$
For each $Z \in \D (\tri) \cap \C_>^{3n},$ $M$ has a (possibly incomplete) hyperbolic structure, such that the topological ideal triangulation is isotopic to a hyperbolic
  ideal triangulation. In particular, $M$ is irreducible and atoroidal, all edges in $M$ are essential, and the structure is complete if and only if the completeness equation $\pee(Z)=1$ is satisfied.
\end{thm}

The situation is more subtle if one considers solutions in $ \D (\tri) \cap \C_\ge^{3n},$ the so-called \emph{partially flat} ideal hyperbolic triangulations. See Petronio and Weeks \cite{PW} for details.


\section{Ideal points and normal surfaces}
\label{sec:Ideal points and normal surfaces}

Since $\D (\tri)$ is a variety in $(\C \setminus\{ 0\})^{3n},$ Bergman's construction in \cite{be} can be used to define its set of ideal points. Let
\begin{align*}
  Z &= (z_1,{\ldots} ,z''_n),\\
  \log |Z| &= (\log |z_1|, \ldots , \log |z''_n|),\\
  u(Z) &= \frac{1}{\sqrt{1 + (\log |z_1|)^2 + {\ldots} + (\log |z''_n|)^2}}.
\end{align*}
The map $\D (\tri) \to B^{3n}$ defined by $Z \to u(Z)\log |Z|$ is continuous, and the \emph{logarithmic limit set} $\D_\infty (\tri)$ is the set of limit points on $S^{3n-1}$ of its image. Thus, for each 
$\xi \in \D_{\infty}(M)$ there is a sequence $\{ Z_i\}$ in $\D (M)$ such that
\begin{equation*}
   \lim_{i\to \infty} u(Z_i) \log |Z_i| = \xi.
\end{equation*}
The sequence $\{ Z_i\}$ is said to \emph{converge} to $\xi,$ written $Z_{i} \to \xi,$
and $\xi$ is called an ideal
point of $\D (\tri).$ Whenever an
edge invariant of a tetrahedron converges to one, the other two edge
invariants ``blow up''. Thus, an ideal point of $\D (\tri)$ is approached
if and only if a tetrahedron degenerates.

Since the Riemann sphere is compact, there is a subsequence, also denoted by $\{ Z_i\},$ with the property that each shape parameter converges in $\C \cup \{ \infty \}.$ In this case $\{ Z_i\}$ is said to \emph{strongly converge} to $\xi.$ If $\xi$ has rational coordinate ratios, a strongly convergent sequence may be chosen on a curve in $\D (\tri)$ according to Lemma 6 in \cite{tillus_ei}.


\subsection{Equivalent descriptions}
\label{log lim}

Let $V$ be a subvariety of $(\C \setminus\{ 0\})^m$ defined by an ideal $J,$ and let $\C [X^{\pm}]=\C [X_1^{\pm 1}, \ldots , X_m^{\pm 1}].$ Bergman \cite{be} noticed that the logarithmic limit set has the two following equivalent descriptions. It is:

1. the set of $m$--tuples $(-v(X_1), \ldots , -v(X_m))$ as $v$ runs over all real--valued valuations on $\C [X^{\pm}]/J$ satisfying $\sum v(X_i)^2 = 1,$ and

2. the intersection over all non--zero elements of $J$ of the spherical duals of their Newton polytopes (see \cite{be, tillus_ei} for details). As noted in \cite{be}, the spherical dual of the convex hull of a finite subset $F \subset \Z^m$ of cardinality $r$ is a finite union of convex spherical polytopes. It is the union over all $\alpha_0, \alpha_1 \in F$ of the set of $\xi$ satisfying the $2r$ inequalities
\begin{equation}\label{eq: com log lim}
  \alpha_0 \cdot \xi \ge \alpha \cdot \xi 
  \qquad\text{and}\qquad
  \alpha_1 \cdot \xi \ge \alpha \cdot \xi,
\end{equation}
where $\alpha$ ranges over $F.$ This will be used in calculations below.


\subsection{The relationship with normal surface theory}

Recall the description of the hyperbolic gluing equation (\ref{eq:glue}) of
$e_j$:
\begin{equation*}
    1 = \prod_{i=1}^n z_i^{a_{ij}} {(z'_i)}^{a'_{ij}} 
{(z''_i)}^{a''_{ij}},
\end{equation*}
and the Q--matching equation \cite{part1}, Section 2.9, of $e_j$:
\begin{equation*}
    0 = \sum_{i=1}^n (a''_{ij} - a'_{ij}) q_i
                     + (a_{ij} - a''_{ij}) q'_i
+ (a'_{ij} - a_{ij}) q''_i.
\end{equation*}
To state the relationship between these sets of equations, let
\begin{equation}
    A = \begin{pmatrix}
           a_{11} & a'_{11} & a''_{11} & a_{21} & {\ldots} & a''_{n1} 
\\
\vdots &        &         &        &          & \vdots   \\
a_{1n} & {\ldots} &       &        &          & a''_{nn} \\
         \end{pmatrix},
\end{equation}
be the exponent matrix of the hyperbolic gluing equations, and
\begin{equation}
    B = \begin{pmatrix}
           a''_{11}-a'_{11} & a_{11}-a''_{11} & a'_{11}-a_{11} & 
{\ldots}
& a'_{n1}-a_{n1} \\
\vdots         &         &        &          & \vdots   \\
a''_{1n}-a'_{1n} & {\ldots}        &        &          &
a'_{nn}-a_{nn} \\
         \end{pmatrix},
\end{equation}
be the coefficient matrix of the Q--matching equations. Let
\begin{equation}
    C_1 = \begin{pmatrix}
       0 & 1 & -1 \\
       -1 & 0 & 1 \\
       1 & -1 & 0
       \end{pmatrix},
\end{equation}
and let $C_n$ be the $(3n \times 3n)$ block diagonal matrix with $n$ copies of
$C_1$ on its diagonal. Then $C_n = - C_n^T,$ and
\begin{equation}
    A C_n = B.
\end{equation}

The deformation variety $\D (\tri)$ is not defined by a principal ideal,
hence its logarithmic limit set $\D_\infty(\tri)$ is in general not
directly determined by its defining equations (see \cite{tillus_ei} for
details).  However, it is contained in the intersection of the spherical duals
of its defining equations:
\begin{equation} \label{comb:tent log lim}
      \D_\infty(\tri) \subseteq \D_{\text{pre-}\infty}(\tri) = 
\overset{n}{\underset{i=1}{\bigcap}}
      \big( Sph(g_i) \cap Sph(p_i) \cap Sph(p'_i) \cap Sph(p''_i) 
\big).
\end{equation}
The set $\D_{\text{pre-}\infty} (\tri)$ is termed a \emph{tropical pre-variety}. In order to give a description of this set, consider first the intersection:
\begin{equation}
    S_n 
      = \overset{n}{\underset{i=1}{\bigcap}}
      \big( Sph(p_i) \cap Sph(p'_i) \cap Sph(p''_i) \big).
\end{equation}
A calculation using the equations (\ref{eq: com log lim}) shows that $S_n$ is the set of all points which are made up of $n$ coordinate triples, each of the form
$(0,x,-x),$ $(-x,0,x)$ or $(x,-x,0),$ where $x \ge 0.$  (If $x>0$ these
correspond to the cases where $z' \to \infty,$ $z'' \to \infty$ or $z \to
\infty$ respectively.) Similarly, one obtains that each hyperbolic gluing
equation gives rise to the intersection of $S_n$ with a hyperplane, leaving
all $\xi \in S_n$ such that
\begin{equation}
    (a_{1j}, a'_{1j}, {\ldots} , a''_{nj})^T \cdot \xi = 0.
\end{equation}
Thus, $\D_{\text{pre-}\infty}(\tri)$ is the intersection of $S_n$ with the nullspace of $A.$


\begin{pro} \label{comb:homeo}
Let $M$ be the interior of an orientable, connected, compact 3--manifold with non-empty boundary consisting of tori, and $\tri$ be an ideal triangulation of $M.$ The set $\D_{\text{pre-}\infty}(\tri)$ is homeomorphic with the projective admissible solution space $\N(\tri)$ of spun-normal surface theory. In particular, $\D_\infty (\tri)$ is homeomorphic with a closed subset of $\N (\tri).$
\end{pro}

\begin{proof}
  The set $\N(\tri)$ is the collection of all elements in the nullspace of
  $B$ with the property that at most one quadrilateral type has non--zero
  coordinate for each tetrahedron, all coordinates are $\ge 0,$ and their sum
  is equal to $1.$ This set may be projected from the unit simplex onto the
  sphere of radius $1/\sqrt{2}$ centered at the origin in $\RR^{3n},$ and, for
  simplicity, this set is also denoted by $\N(\tri).$
  
  The map $C_n^T$ takes $\N(\tri)$ to the unit sphere $S^{3n-1},$ where
  $\D_{\text{pre-}\infty}(\tri)$ is found, since the following correspondence between the
  $i$--th coordinate triples holds:
  \begin{equation}\label{comb:C}
    \begin{pmatrix} 0 \\ x \\ -x \end{pmatrix}
    = C^T
    \begin{pmatrix} x \\ 0 \\ 0 \end{pmatrix}
                  \qquad
    \begin{pmatrix} -x \\ 0 \\ x \end{pmatrix}
    = C^T
    \begin{pmatrix} 0 \\ x \\ 0 \end{pmatrix}
                \qquad
    \begin{pmatrix} x \\ -x \\ 0 \end{pmatrix}
    = C^T
    \begin{pmatrix} 0 \\ 0 \\ x \end{pmatrix}.
  \end{equation}
  Thus, if $N \in \N(\tri)$ satisfies $|| N ||^2 = 1/2,$ then $|| C_n^T N
  ||^2 = 2 ||N||^2 = 1.$  Furthermore, given $N \in \N(\tri),$ one has
  $C_n^T N \in S_n,$ the set containing $\D_{\text{pre-}\infty}(\tri).$ Now $0 = B N = A
  (C_nN) = -A(C_n^TN)$ implies that $C_n^T N \in \D_{\text{pre-}\infty}(\tri).$ Thus, there
  is a linear map $\N(\tri) \to \D_{\text{pre-}\infty}(\tri).$
  
  The kernel of $C_n^T$ is generated by the vectors with $(1,1,1)^T$ in the
  $i$--th triple and $0$ in the other positions. It follows that different
  admissible solutions cannot differ by an element in the kernel. The linear
  map is therefore 1--1.
  
  Since every element in $S_n$ has a unique inverse image under $C_n^T,$ any
  $\xi \in \D_{\text{pre-}\infty}(M)$ can be taken to a normal $Q$--coordinate $N(\xi)$ using
  (\ref{comb:C}). Thus, $\xi = C_n^T N(\xi),$ and hence $BN(\xi)=0.$  This
  shows that the map is onto, and in fact, that there is a well--defined
  inverse mapping.
  
  The last claim follows since $\D_\infty (\tri)$ is a closed subset of $\D_{\text{pre-}\infty}(\tri).$
\end{proof}

The proof of Proposition \ref{comb:homeo} shows that $\D_\infty (\tri)$ is homeomorphic to a closed subset of $\N (\tri)$ in a canonical way. For $\xi \in \D_\infty (\tri),$ let $N(\xi)$ be the unique normal $Q$--coordinate such that $\xi = C_n^T N(\xi).$ If $\xi$ has rational coordinate ratios, one can associate a unique spun-normal surface $S(\xi)$ to it as follows. By assumption, there is $r>0$ such that $rN(\xi)$ is an integer solution, and hence corresponds to a unique spun-normal surface without vertex linking components (see Theorem 2.4 in \cite{part1}). One then requires $r$ to be minimal with respect to the condition that the surface is 2--sided. The properties studied below are independent of the choice of $S(\xi),$ but it will be convenient to refer to a surface.


\section{Boundary curves of essential surfaces}
\label{Holonomy variety and essential surfaces}

In this section, the derivative of the holonomy of \cite{nz} is related to the linear functional $\nu$ of \cite{part1}, and used to show that the boundary curves of the spun-normal surfaces in $\D_\infty(\tri)$ are the boundary curves of essential surfaces which are strongly detected by the character variety. 
Recall that one calls the boundary slope of an essential surface associated to an ideal point of the character variety \emph{strongly detected} if no closed surface can be associated to that particular ideal point.


\subsection{Essential surfaces and boundary curves}

An \emph{(embedded) surface} $S$ in the topologically finite 3--manifold $M = \text{int} (\overline{M})$ will always mean a 2--dimensional PL submanifold of $M$ with the property that its closure $\overline{S}$ in $\overline{M}$ is \emph{properly embedded} in $\overline{M},$ that is, a closed subset of $\overline{M}$ with $\partial \overline{S} = \overline{S} \cap \partial \overline{M}.$ A surface $S$ in $M$ is said to be \emph{essential} if its closure $\overline{S}$ is essential in $\overline{M}$ as described in the following definition:
\begin{defn*} \cite{sh1}
A surface $\overline{S}$ in a compact, irreducible, orientable 3--manifold $\overline{M}$ is said to be essential if it has the following five properties:
    \begin{enumerate}
    \item $\overline{S}$ is bicollared;
    \item the inclusion $\pi_1(\overline{S}_i) \to \pi_1(\overline{M})$ is injective for every
      component $\overline{S}_i$ of $\overline{S}$;
    \item no component of $\overline{S}$ is a 2--sphere;
    \item no component of $\overline{S}$ is boundary parallel;
    \item $\overline{S}$ is nonempty.
    \end{enumerate}
\end{defn*}
The boundary curves of $\overline{S},$ $\overline{S} \cap \partial\overline{M},$ are also called the \emph{boundary curves} of $S.$


\subsection{Simplicial version of $\mathbf{\nu(\gamma)}$}

\begin{figure}[t]
\psfrag{z}{{\small $z$}}
\psfrag{z'}{{\small $z'$}}
\psfrag{z"}{{\small $z''$}}
\psfrag{q}{{\small $q$}}
\psfrag{q'}{{\small $q'$}}
\psfrag{q"}{{\small $q''$}}
\psfrag{g}{{\small $\gamma$}}
\psfrag{h}{{\small $\gamma'$}}
\begin{center}
      \subfigure[$\mu$ and $\nu$]{
        \includegraphics[width=4cm]{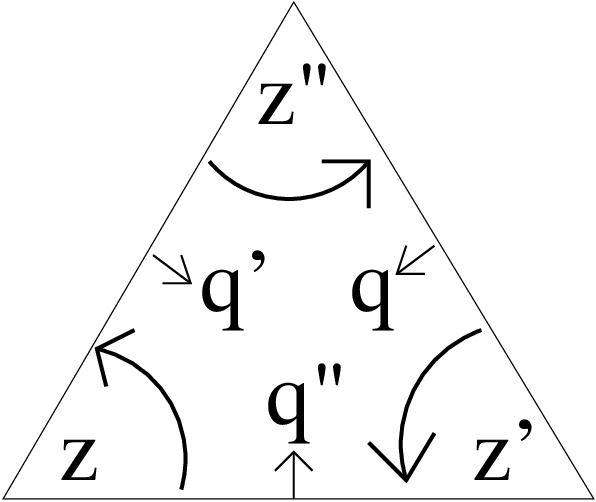}}
       \qquad \qquad
     \subfigure[$\nu^{simp}(\gamma) = \nu(\gamma')$]{
         \includegraphics[height=4cm]{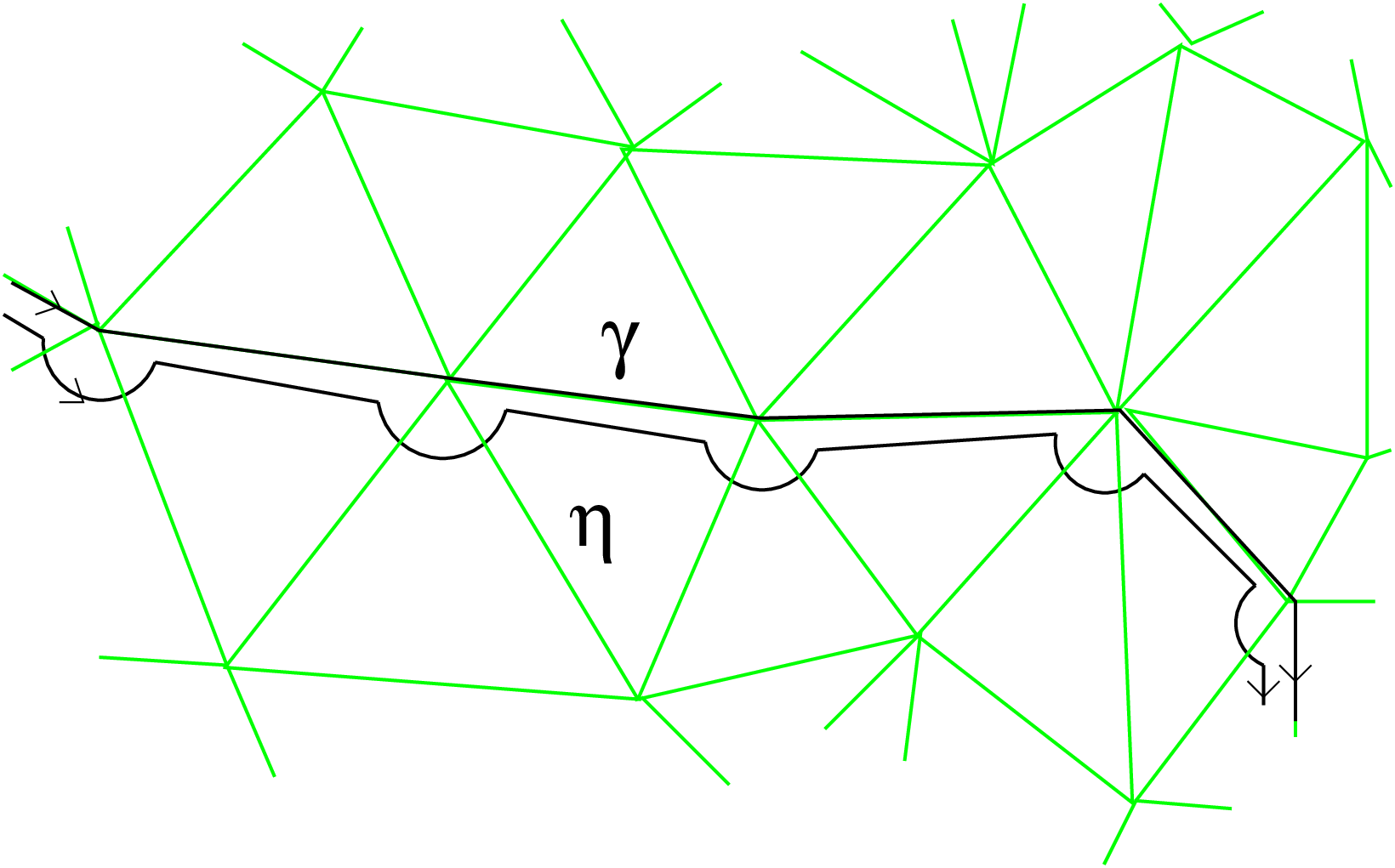}}
\end{center}
    \caption{$\mu$ and $\nu$}
    \label{fig:mu and nu}
\end{figure}

It will be convenient to have a simplicial definition of the functional $\nu(\gamma)$ defined in \cite{part1}, Section 3.1. Let $\Delta$ be a triangle in the induced triangulation $\tri_i$ of $T_i.$  If $v,u,t$ are the vertices of $\Delta$ in clockwise ordering (as viewed from the cusp), and $q_0$ is the $Q$--modulus of $[v,t],$ and $q_1$ is the $Q$--modulus of $[v,u],$ then define the \emph{$Q$--modulus} of $v$ to be $q_0-q_1$ (with respect to $\Delta$). E.g. if $\Delta$ is the triangle of Figure \ref{fig:mu and nu}(a), then the $Q$--modulus of the vertex with label $z$ is $q''-q'$ with respect to $\Delta.$ Note that for a fixed vertex, the sum of moduli with respect to all triangles containing it is equal to the corresponding $Q$--matching equation.

If $\gamma$ is an oriented simplicial path, let $\nu^{simp}(\gamma)$ be the sum of $Q$--moduli of vertices of triangles touching $\gamma$ to the right. If $\gamma'$ is the right hand boundary component of a small regular neighbourhood of $\gamma,$ oriented in the same way as $\gamma,$ then $\nu(\gamma') = \nu^{simp}(\gamma)$; hence put $\nu^{simp}=\nu.$  Let $N$ be a solution to the $Q$--matching equations, and $\nu_N(\gamma )$ be the evaluation of $\nu(\gamma )$ at $N.$

\begin{lem} [\cite{part1}, Lemma 3.1] \label{boundary:homo}
  Let $N$ be a solution to the $Q$--matching equations.  The number
  $\nu_N(\gamma ) \in \RR$ depends only on the homotopy class of $\gamma$ and
  defines a homomorphism $\nu_N : \pi_1(T_i) \to (\RR , +).$\qed
\end{lem}

Let $\overrightarrow{\nu(\gamma)}$ be the coefficient (row) vector of $\nu(\gamma),$ and $\overrightarrow{\mu(\gamma)}$ be the exponent (row) vector of $\mu(\gamma),$ where $(q_1,q'_1,q''_1,{\ldots} ,q''_n)$ and $(z_1,z'_1,z''_1,{\ldots} ,z''_n)$ are the respective coordinate systems. Then:

\begin{lem} \label{lem:nu&mu}
$\overrightarrow{\nu(\gamma)} = \overrightarrow{\mu(\gamma)} C_n^T$
\end{lem}

\begin{proof}
  It it sufficient to verify the relationship for vertex moduli, and hence for
  the vertex labels of the triangle in Figure \ref{fig:mu and nu}. If $\gamma$
  touches the vertex with label $z$ to the right, then the contribution to
  $\mu(\gamma)$ is $z,$ and the contribution to $\nu(\gamma)$ is $q''-q'.$
  Similarly, $z'$ corresponds to $q-q''$ and $z''$ to $q'-q.$
  Describing this relationship for coordinate triples using $C_1$ gives:
  $$(1,0,0)C_1^T = (0,-1,1), \quad (0,1,0) C_1^T = (1,0,-1),\quad (0,0,1)C_1^T =
  (-1,1,0).$$ 
  This proves the lemma.
\end{proof}


\subsection{Strongly detected boundary curves}

The following proposition only concerns the boundary curves of a spun-normal surface, and its proof \emph{uses} the known results from \cite{cs}. Section \ref{The geometric construction} gives an independent proof of the fact that the non-compact spun-normal surfaces involved are indeed non-trivial.

\begin{pro}\label{comb:essential prop}
Let $M$ be the interior of a compact, connected, irreducible, orientable 3--manifold with non-empty boundary consisting of tori, and let $\tri$ be an ideal triangulation of $M$ with the property that all edges are essential. Let $\xi \in \D_\infty(\tri)$ be a point with rational coordinate ratios, and assume that $S(\xi)$ is not closed. Then there is an essential surface in $M$ dual to an ideal point of a curve in the character variety of $M$ which has (up to projectivisation) the same boundary curves as $S(\xi).$ Moreover, these boundary curves are strongly detected.
\end{pro}

\begin{proof}
  Fix a basis $\{\m_i, \l_i\}$ for each boundary torus $T_i$ of $M,$ and
  denote the resulting coordinates for $\PEi (M)$ by $(M_1,L_1,{\ldots} ,M_h,
  L_h).$ Let $\xi \in \D_\infty(\tri)$ be an ideal point with rational
  coordinate ratios. Lemma 6 of \cite{tillus_ei} provides a curve $C$ in
  $\D(\tri)$ such that $\xi \in C_\infty.$  Moreover, there is $\alpha >0$
  such that $\alpha \xi$ defines a normalised, discrete, rank 1 valuation
  $v_\xi$ on $\C (C)$:
  \begin{equation*}
    \alpha \xi = (-v_\xi(z_1),{\ldots} ,-v_\xi(z''_n)).
  \end{equation*}
  Indeed, one may choose $\alpha>0$ such that the above is an integer valued
  vector whose entries have no common divisor. Then the normal $Q$--coordinate
  of $S(\xi)$ is $\alpha N(\xi )$ or $2\alpha N(\xi ).$ Therefore denote $\alpha \xi$ and $\alpha
  N(\xi )$ by $\xi$ and $N(\xi )$ respectively.
  
  For each boundary torus $T,$ and each $\gamma \in \im (\pi_1(T) \to
  \pi_1(M)),$ there are $g_1,{\ldots} ,g''_n \in \Z$ such that:
  \begin{equation*}
    \mu(\gamma ) =
       \prod_{i=1}^n z_i^{g_{i}} {(z'_i)}^{g'_{i}} {(z''_i)}^{g''_{i}},
     \text{ and }
    \overrightarrow{\mu(\gamma)} = (g_1,{\ldots} ,g''_n) \in \Z^{3n}.
  \end{equation*}
  Recall that $\xi = C_n^T N(\xi)$ and $C_n = - C_n^T.$ Using this and Lemma
  \ref{lem:nu&mu}, one has:
  \begin{align*}
  v_\xi(\mu (\gamma) ) &= v_\xi \bigg( \prod_{i=1}^n z_i^{g_{i}}
  {(z'_i)}^{g'_{i}} {(z''_i)}^{g''_{i}} \bigg) = \sum_{i=1}^n g_i v_\xi(z_i) +
  g'_iv_\xi(z'_i) +
  g''_iv_\xi(z''_i)\\
  &= - \overrightarrow{\mu(\gamma)} \cdot \xi = - \overrightarrow{\mu(\gamma)}
  \cdot C_n^T N(\xi)
  = - \overrightarrow{\nu(\gamma)} \cdot N(\xi)\\
  &= - \nu_{N(\xi)}(\gamma) = \nu_{N(\xi)}(\gamma^{-1}).
  \end{align*}
  Thus, restricted to the boundary, $v_\xi\mu = -\nu_{N(\xi)},$ and in
  particular, since $S(\xi)$ is not closed, the eigenvalue of at least one
  peripheral element blows up. This implies that the restriction $\pee : C \to
  \PEi (M)$ is not constant, and its image is therefore a curve $C' \subset
  \PEi (M).$ Denote the ideal point of $C'$ corresponding to $\xi$ by $\xi'.$
  A corresponding normalised, discrete, rank 1 valuation $v'$ of $\C (C')$ at
  $\xi'$ is obtained as follows. If
  \begin{equation}\label{comb:val}
    (v_\xi \mu(\m_1), v_\xi\mu(\l_1),{\ldots} ,v_\xi \mu(\m_h ),
    v_\xi\mu(\l_h)) 
  \end{equation}
  contains a pair of coprime integers, define $v'(M_i) = v_\xi \mu (\m_i),$
  and $v'(L_i) = v_\xi \mu (\l_i).$  Otherwise, let $d$ denote the greatest
  common divisor of the entries in (\ref{comb:val}), and then define $v'(M_i)
  = \frac{1}{d}v_\xi \mu (\m_i),$ and $v'(L_i) = \frac{1}{d}v_\xi \mu (\l_i).$
  In either case, one obtains a valuation of $\C (C')$ with the desired
  properties.
  
  Culler-Shalen theory can be applied using $\xi'$ and $v'.$ Lemma 14 of
  \cite{tillus_ei} yields that the projectivised boundary curves of an
  essential dual surface are given by:
  \begin{align*}
  & [v' (\l_1), -v'(\m_1),{\ldots} ,v' (\l_h ), -v'(\m_h)]\\
  = &[v_\xi (\mu\l_1),-v_\xi(\mu\m_1),{\ldots} ,v_\xi (\mu\l_h ),
  -v_\xi(\mu\m_h)]\\
  = & [-v_{N(\xi)} (\l_1), v_{N(\xi)}(\m_1),{\ldots} ,-v_{N(\xi)} (\l_h ),
  v_{N(\xi)}(\m_h)].
  \end{align*}
  This proves the claim since the latter gives the projectivised boundary
  curves of $S(\xi)$ according to equation (4.1) in \cite{part1}.
\end{proof}


\section{The leaf space}
\label{sec:The dual tree}


Throughout this section, let $M$ be the interior of a compact, connected 3--manifold with non-empty boundary, and $\tri$ be an ideal triangulation of $M.$ Every element, $N,$ of the projective admissible solution space $\N(\tri)$ of \cite{part1} is interpreted as a transversely measured singular codimension--one foliation of $M.$ This foliation is lifted to the universal cover $\widetilde{M}$ and the leaf space is turned into an $\RR$--tree, $\tree_N,$ on which the fundamental group of $M$ acts by isometries (Theorem~\ref{thm:general tree}).


\subsection{Singular foliations}

\begin{figure}[t]
\begin{center}
      \subfigure[The singular codimension-one foliation]{\label{fig:foliation}
        \includegraphics[width=12cm]{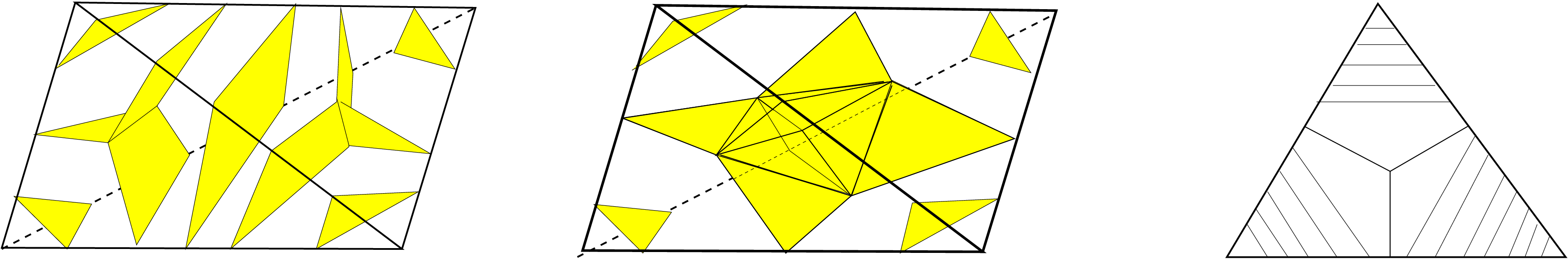}}
      \\
           \subfigure[The dual spine]{\label{fig:spine}
         \includegraphics[width=12cm]{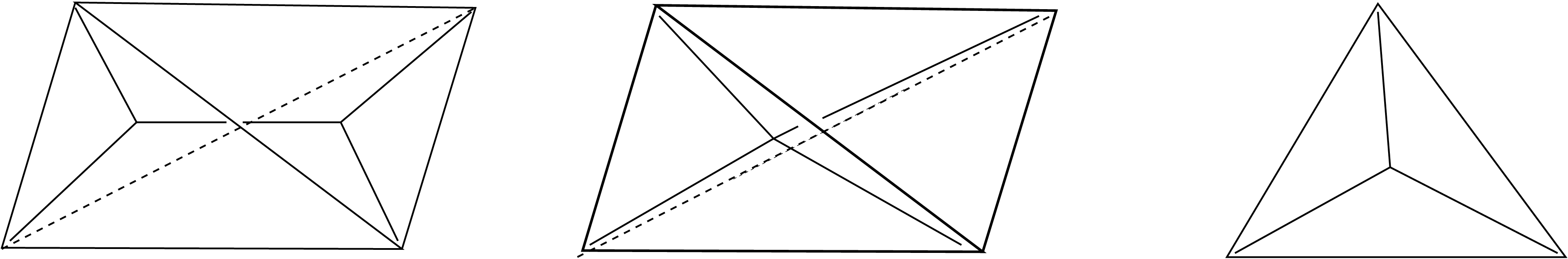}}
\end{center}  \caption{The transversely measured singular foliation}
\end{figure}

The element $N \in \N (\tri)$ determines a singular foliation of each ideal tetrahedron, $\sigma,$ as follows. If all quadrilateral coordinates of $N$ supported by $\sigma$ are zero, then the foliation consists of triangles and a single singular leaf \emph{(butterfly).} Otherwise there is a unique non-zero quadrilateral coordinate, and the foliation consists of quadrilaterals of that type, triangles and two singular leaves \emph{(wings)}. The two kinds of singular foliation are indicated in Figure \ref{fig:foliation}. The pattern on each face is topologically the same, and determines a singular foliation thereof. There also is a natural transverse measure, obtained as follows. Consider a \emph{dual spine} in $\sigma$ which is transverse to the foliation. It naturally inherits the structure of a metric tree: it consists of two copies of $[0, \infty)$ attached at each end point of a (possibly degenerate) closed interval $[0,k],$ where $k$ is the maximum over all quadrilateral coordinates of $N$ supported by $\sigma.$ The metric tree then defines a unique transverse measure to the foliation.

There is a unique way to identify the foliations of adjacent tetrahedra such that the leaves and transverse measures match up across the faces. The $Q$--matching equations ensure that the leaves close up around the edges of the triangulation, so the singular foliations of the tetrahedra glue up to give a singular foliation of $M$ with finitely many singular leaves. If all coordinate ratios in $N$ are rational, then Theorem 2.4 of \cite{part1} implies that all leaves are closed and proper, giving a foliation by spun-normal surfaces together with finitely many singular leaves.

If all ends of $M$ are tori, the homomorphism $\nu_N$ determines the behaviour of the singular foliation of
$M$ near the cusps as in \cite{part1}, Section 4.1.
Thus, if $\nu_N(\l_i)=\nu_N(\m_i)=0$ for some $i,$ then the foliation near the
$i$--th cusp of $M$ is topologically of the form $T^{2}\times (0, 1).$
Otherwise, the fraction $-\nu_N(\l_i)/\nu_N(\m_i)$ determines the slope of the
intersection of the leaves with the cusp, and hence the topology of the
induced foliation of $\partial M.$


\subsection{Equivalence relations}
\label{sec:Equivalence relation}

Let $N$ be an admissible solution to the Q--matching equations, and denote by $\fol$ the transversely measured singular codimension--one foliation of $M$ defined by $N.$ In each tetrahedron place the corresponding dual spine. The universal cover $\widetilde{M}$ is given the ideal triangulation induced by $\tri,$ so that the covering map $p\co \widetilde{M} \to M$ is simplicial and $\pi_1(M)$ acts simplicially on $\widetilde{M}.$ Then $\fol$ lifts to a transversely measured singular codimension--one foliation $\widetilde{\fol}$ of $\widetilde{M}$ which is invariant under the group action. 

\begin{figure}[t]
\psfrag{p}{{\small $p$}}
\psfrag{q}{{\small $q$}}
\psfrag{r}{{\small $r$}}
\psfrag{s}{{\small $s$}}
\psfrag{R}{{\small $r=\max\{p,q\}-\min\{p,q\}$}}
\psfrag{S}{{\small $s=\min\{p,q\}$}}
  \begin{center}
    \subfigure[Case 1]{
      \includegraphics[height=5cm]{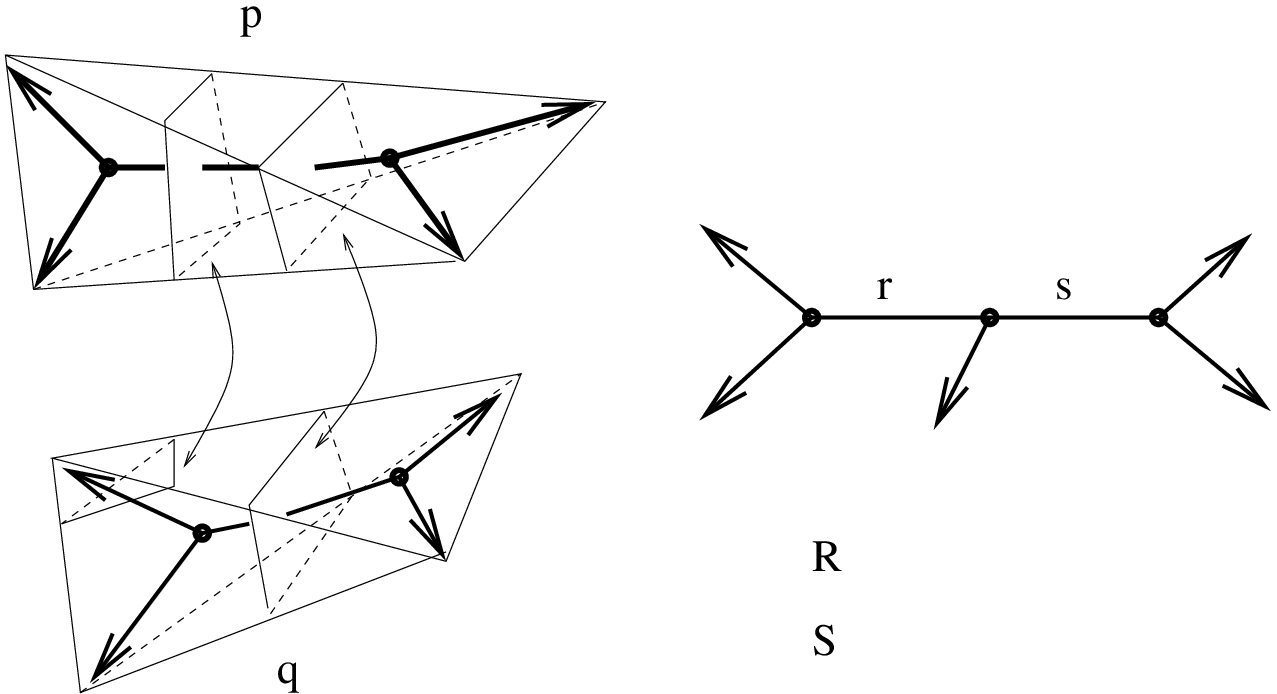}}
      \\
    \subfigure[Case 2]{
      \includegraphics[height=5cm]{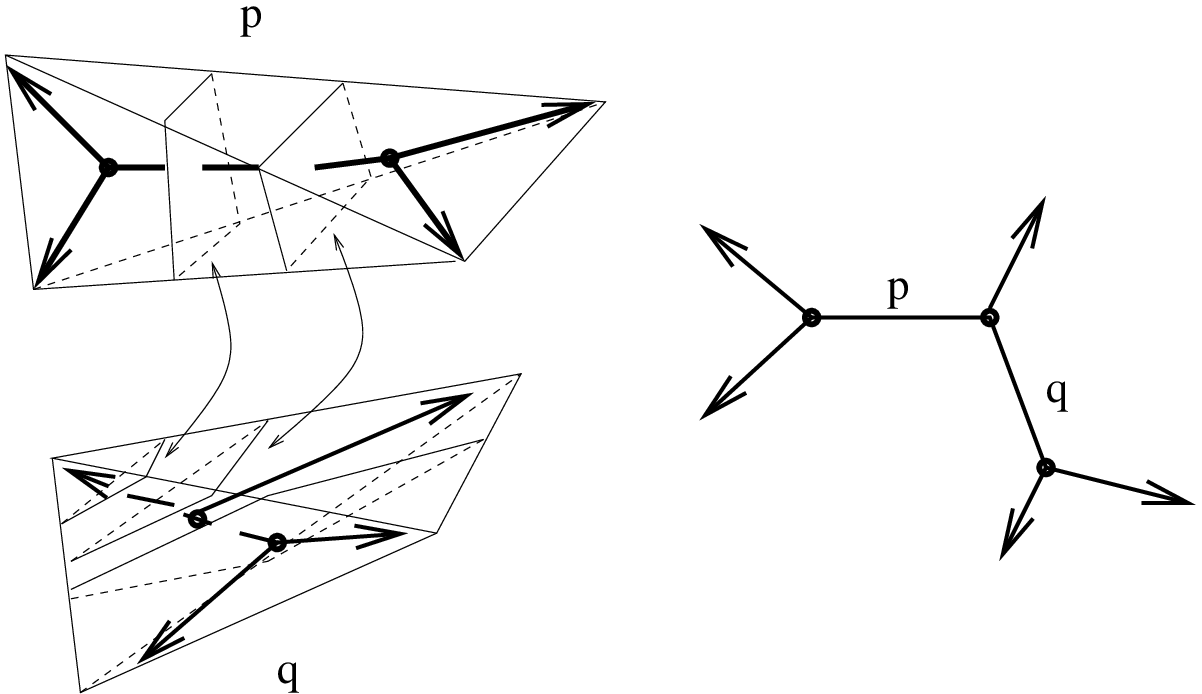}}
     \caption{Identification of dual spines with $p, q \in [0, \infty)$}
   \label{fig:dual_cases}
  \end{center}
\end{figure}

The leaf spaces $M/\fol$ and $\widetilde{M}/\widetilde{\fol}$ can be obtained by introducing an equivalence relation on the disjoint union of all dual spines, which is generated as follows. Each point on a dual spine corresponds to a leaf, and the identification of two ideal 3--simplices along a common face induces identifications of their dual spines according to how the leaves and the transverse measures match up across the face, as shown in Figure \ref{fig:dual_cases}. Indeed, the identification of points on adjacent dual spines coming from a face is uniquely determined by the way the ideal vertices match up: they only occur along the spines minus the interior of the intervals $[0,\infty)$ corresponding to the ideal vertices not belonging to the face. The results are of the shape of a $Y,$ each of whose ends corresponds to an ideal vertex of the face. Across the face, the singularities of the corresponding $Y$'s match up, and the ends are identified according to the face pairing and the transverse measure.

The foliation induces a pseudo--metric on $\widetilde{M},$ which descends to a pseudo--metric on $\widetilde{M}/\widetilde{\fol}.$ Any path $\gamma$ in $\widetilde{M}$ inherits a measure $\mu_N (\gamma),$ which is the total mass of the transverse measure along it. The measure of a path contained in a leaf is zero. Call a path $\gamma$ \emph{admissible} if it is made up of pieces each of which is either contained in a leaf of $\widetilde{\fol}$ or transverse to $\widetilde{\fol}.$ Any path is homotopic to an admissible path by a homotopy which does not increase the measure. For $x,y \in \widetilde{M},$ define $d_N(x,y)= \inf \mu_N(\gamma),$ where the infimum is taken over all admissible paths from $x$ to $y.$ The pseudo--metric $d_N$ on $\widetilde{M}$ descends to a pseudo--metric $d_N'$ on $\widetilde{M}/\widetilde{\fol}.$ Two points $x,y \in \widetilde{M}$ are equivalent, written $x \sim y,$ if $d_N(x,y)=0.$ Denote by $\tree_N$ the space $\widetilde{M}/\sim$ with the induced metric. Then $\widetilde{M}/\widetilde{\fol}$ is isomorphic to $\tree_N$ if and only if $d_N'$ is a metric, and this is equivalent to $\widetilde{M}/\widetilde{\fol}$ being Hausdorff. 

Since the foliation and the transverse measure are invariant under the action of $\pi_1(M)$ by deck transformations on $\widetilde{M},$ there is an induced isometric action of $\pi_1(M)$ on $\tree_N.$ This is the unique action which makes the quotient map $\widetilde{M}\to \tree_N$ $\pi_1(M)$--equivariant and continuous.

\begin{thm}\label{thm:general tree}
Let $M$ be the interior of a compact, connected 3--manifold with non-empty boundary, and $\tri$ be an ideal triangulation of $M.$ For each admissible solution $N$ of the $Q$--matching equations, $\tree_N$ is an $\RR$--tree on which $\pi_1(M)$ acts by isometries.
\end{thm}

\begin{proof}
It remains to show that $\tree_N$ is an $\RR$--tree. The argument given is due to Paulin \cite{P}. According to \cite{ab} it is necessary and sufficient to show that $\tree_N$ is connected and $0$--hyperbolic. Let $x,y \in \tri_N.$ Then there are points $\tilde{x}, \tilde{y} \in \widetilde{M}$ corresponding to $x$ and $y.$ Since $\widetilde{M}$ is connected, there is an admissible path from $\tilde{x}$ to $\tilde{y}.$ This path descends to a path in $\tree_N,$ showing that it is connected.

Now $\tree_N$ is $0$--hyperbolic if for all $x_0,x_1,x_2,x_3 \in \tree_N,$ one has
\begin{equation}\label{eq:0-hyp}
d(x_0,x_2) + d(x_1,x_3) \le \max \{ d(x_0,x_1) + d(x_2,x_3), d(x_0,x_3) + d(x_1,x_2)\}.
\end{equation}
Regard indices in the following as integers modulo four. Let $x_i \in \tree_N,$ and choose preimages $\tilde{x}_i \in \widetilde{M}.$ By definition, given $\varepsilon > 0,$ there are admissible paths $\tilde{\alpha}_i$ connecting $\tilde{x}_i$ to $\tilde{x}_{i+1}$ with the property that
\begin{equation*}
\mu_N(\tilde{\alpha}_i)\le d_N(\tilde{x}_i, \tilde{x}_{i+1}) + \varepsilon.
\end{equation*}
Since $\widetilde{M}$ is simply connected, there is a PL map $f$ from the disc $D$ to $\widetilde{M}$ such that the preimages $\overline{x}_i$ of $\tilde{x}_i$ are ordered according to indices on the boundary of $D$ and such that the subarc $\overline{\alpha}_i$ on $\partial D$ connecting $\overline{x}_i$ to $\overline{x}_{i+1}$ is mapped injectively onto $\tilde{\alpha}_i.$ 

Without loss of generality, it may be assumed that $f$ is in general position, so that $\fol$ pulls back to a transversely measured singular foliation $\overline{\fol}$ of $D.$ The disc $D$ will be viewed as a square with sides $\overline{\alpha}_i;$ the inequality (\ref{eq:0-hyp}) compares the measures of opposite sides to the measure on the diagonals. The Poincar\'e recurrence theorem assures that there are at most finitely many singularities in $\overline{\fol}.$ A leaf transverse to $\partial D$ cannot accumulate in the interior of $D,$ and the measure of all points on $\overline{\alpha}_i$ which are contained on leaves having two endpoints on $\overline{\alpha}_i$ is bounded above by $2 \varepsilon.$ Thus, $\overline{\alpha}_i$ can be divided into three subarcs with the property that (apart from the set of leaves going back to $\overline{\alpha}_i$) each leaf meeting the first (second, third) subarc has an endpoint on $\overline{\alpha}_{i+1}$ ($\overline{\alpha}_{i+2},$ $\overline{\alpha}_{i+3}$). Note that if the measure of leaves connecting $\overline{\alpha}_i$ to $\overline{\alpha}_{i+2}$ is non--zero, then the measure of leaves connecting 
$\overline{\alpha}_{i+1}$ to $\overline{\alpha}_{i+3}$ is zero. The inequality (\ref{eq:0-hyp}) now follows, since $\varepsilon$ was arbitrary.
\end{proof}


\subsection{Dual tree of a surface}
\label{subsec:dual tree}

If $S$ is a 2--sided surface in $M,$ then there is a well--defined \emph{dual simplicial graph} $\mathcal{G}_S$ associated with $S,$ whose edges and vertices are in bijective correspondence with the components of $S$ and $M-S$ respectively (see \cite{sh1}, Section 1.4). There are a retraction $r\co M \to \mathcal{G}_S$ and an inclusion $i \co \mathcal{G}_S \to M,$ implying that $\mathcal{G}_S$ is connected and that $\pi_1(\mathcal{G}_S)$ is isomorphic to a subgroup and a quotient of $\pi_1(M).$ Let $\tilde{S} = p^{-1}(S)$ and denote the dual graph associated with $\tilde{S}$ by $\tree_S.$ It follows that $\pi_1(\tree_S)$ is trivial, so $\tree_S$ is a simplicial tree. Moreover, the action of $\pi_1(M)$ on $\widetilde{M}$ induces a simplicial action without inversions on $\tree_N;$ this is the unique action which makes the map $\widetilde{M} \to \tree_N$ $\pi_{1}(M)$--equivariant.

\begin{pro}
Let $M$ be the interior of a compact, connected 3--manifold with non-empty boundary, and $\tri$ be an ideal triangulation of $M.$ Let $N$ be an admissible integer solution of the $Q$--matching equations with the property that the (unique) associated spun-normal surface $S$ is 2--sided. Then $\tree_N$ is a simplicial tree on which $\pi_1(M)$ acts simplicially without inversions. Moreover, it is isomorphic to $\widetilde{M}/\widetilde{\fol}$ and contains the dual graph of $\tilde{S}=p^{-1}(S)$ as a simplicial subtree.
\end{pro}

\begin{proof}
Replace the surface $S$ by a surface $S'$ which may have infinitely many boundary parallel components such that it contains infinitely many triangle discs of each type. The dual spine determined by the corresponding solution $N$ of the $Q$--matching equations inherits the structure of a simplicial graph: it consists of two copies of $[0, \infty)$ attached at each end point of a (possibly degenerate) closed interval $[0,k],$ where $k$ is the number of normal quadrilaterals of $S$ in $\sigma.$ Assume that the normal triangles meet the half--open intervals in precisely the half--integer places, and that the normal quadrilaterals meet the closed interval in precisely the half--integer places. The integer places are referred to as the vertices of the spine, and each interval between two vertices as an edge.

The dual graph $\mathcal{G}_{S'}$ is the dual graph $\mathcal{G}_S$ together with a copy of $[0, \infty)$ attached to each vertex corresponding to a component of $M-S$ containing an ideal vertex of $\tri.$ The dual tree of $\tilde{S}=p^{-1}(S)$ is then a subtree of the one of $\tilde{S'}=p^{-1}(S').$

The normal discs divide each tetrahedron into regions, one for each vertex on its dual spine. The components of $\widetilde{M}-\tilde{S'}$ are a partition of $\widetilde{M},$ and there is a corresponding (unique) equivalence relation on the set of regions giving this partition; it is generated by identifying two regions if they meet along a face of a tetrahedron. Similar considerations hold for the spun-normal surface and its cell decomposition by normal discs. Thus, if vertices on spines are identified whenever the corresponding regions are glued to each other, and if edges on spines are identified whenever the corresponding normal discs are glued to each other, then the dual graph $\tree_{S'}$ to $\tilde{S'}$ is obtained. This is exactly the equivalence relation described in Section \ref{sec:Equivalence relation}; hence $\tree_{S'}$ is equivariantly isomorphic to $\widetilde{M}/\widetilde{\fol},$ and this implies that $\widetilde{M}/\widetilde{\fol}$ is isomorphic to $\tree_N.$
\end{proof}

The following consequence is obtained by rescaling:

\begin{cor}\label{cor:trees}
Let $M$ be the interior of a compact, connected 3--manifold with non-empty boundary, and $\tri$ be an ideal triangulation of $M.$ Let $N$ be an admissible solution of the $Q$--matching equations with rational coordinate ratios. Then $\tree_N$ is a simplicial tree with rational edge lengths, and isomorphic to the leaf space $\widetilde{M}/\widetilde{\fol}.$
\end{cor}


\section{Action on the limiting tree}
\label{The geometric construction}

The main results stated in the introduction are proved in this section. The statements of Theorem~\ref{thm:main 1} are contained in Corollaries~\ref{pro: non-trivial action via peripheral} and~\ref{cor:thm 1 part 2}; Part (1) of Theorem~\ref{thm: morgan-shalen} is in Corollary~\ref{cor: non-trivial action 2}, and part (2) of Theorem~\ref{thm: morgan-shalen} is in Proposition~\ref{pro:dual surfaces}.


\subsection{Outline and definitions}
\label{sec:Outline and definitions}

Let $\xi \in \D_{\infty}(\tri),$ and fix a sequence $\{ Z_i\}$ in $\D (\tri)$ which strongly converges to $\xi,$ so 
$$\lim_{i\to \infty} u(Z_i)\log | Z_i| = \xi$$ 
in the notation of Section \ref{sec:Ideal points and normal surfaces}, and each shape parameter converges in $\C \cup \{ \infty \}.$ Denote by $N=N(\xi) \in \N (\tri)$ the image under the natural homeomorphism. If $\xi$ has rational coordinate ratios, denote by $S = S(\xi)$ the associated 2--sided spun-normal surface without boundary parallel components. For each $Z \in \D(\tri),$ there is a map $D_{Z}: \widetilde{M} \to \H^{3}$ as in Section \ref{sec: developing}.  Let $D_{i}$ denote the map corresponding to $Z_{i},$ and $\rho_i : \pi_1(M)\to \PSL$ be the corresponding representation.

The ideal triangulation of $M = \{ \sigma_1,...,\sigma_n\} / \sim$ is lifted equivariantly to an ideal triangulation of the universal cover, so that 
\begin{equation*}
\widetilde{M} =  \{ \gamma\tilde{\sigma}_1,...,\gamma\tilde{\sigma_n} | \gamma \in \pi_1(M)\} / \Psi,
\end{equation*}
where $\Psi$ denotes the face pairing scheme. As in \cite{y}, let $\widetilde{M}_{i}$ be the topological space which is obtained from $\widetilde{M}$ by imbuing each 3--simplex $\sigma$ of $\widetilde{M}$ with the hyperbolic structure determined by $Z_i$ such that all faces are identified isometrically. Thus, if $\sigma(Z_i)$ is an ideal hyperbolic tetrahedron isometric to $D_i(\sigma),$ then taking the disjoint union $\{ \sigma(Z_i) : \sigma \subset \widetilde{M} \}$ together with the isometric face pairings induced from $\Psi$ gives:
\begin{equation}\label{eq: degeneration}
\widetilde{M}_{i} = \{ \sigma(Z_i) : \sigma \subset \widetilde{M} \} / \Psi_i.
\end{equation}
Note that $\widetilde{M}_{i}$ may not be separable, and that there is a $\pi_1(M)$--equivariant proper continuous map $f_i : \widetilde{M}\to \widetilde{M}_i.$

It is shown in Section \ref{GH-limit of tet} that the Gromov-Hausdorff limit of $(\sigma(Z_i), u(Z_i)d)$ is a so--called dual spine $S(\sigma),$ where $d$ denotes the hyperbolic metric. The leaf space of the singular foliation $\fol$ associated to $\xi$ is described in the previous section as the set of dual spines modulo an equivalence relation:
\begin{equation}\label{eq: leaf space}
\widetilde{M}/\widetilde{\fol}= \{ S(\sigma) : \sigma \subset \widetilde{M} \} / \sim.
\end{equation}
This leads to an interpretation of $\widetilde{M}/\widetilde{\fol}$ as a limit of the sequence $(\widetilde{M}_{i}, u(Z_i)d).$

The leaf space may not be Hausdorff; identifying any two non--separable points turns it into an $\RR$--tree $\tree_N.$ There is a map $\widetilde{M}_{i} \to D_{i}(\widetilde{M}) \subseteq \H^3,$ which allows the comparison between the action of $\pi_1(M)$ on $\tree_N$ and the limiting behaviour of the sequence of actions $\rho_i$ on $(\H^3, u(Z_i) d).$ 


\subsection{Geodesic spines}

The image of each ideal (topological) tetrahedron in $\widetilde{M}$ under $D_{i}$ is an ideal hyperbolic tetrahedron with shape parameters in $\C-\{0,1\}$ for each $i.$ To simplify notation, assume that the simplices and parameter triples are ordered such that as $Z_{i}\to\xi,$ the simplices $\sigma_{1},...,\sigma_{k}$ degenerate with $z_{j}\to 1,$ and $\sigma_{k+1},...,\sigma_{n}$ do not degenerate. A geodesic dual spine $S_i(\sigma)$ is constructed in $\sigma(Z_i)$ for each $\sigma \subset \widetilde{M}$ as follows.

\begin{figure}[t]
\psfrag{b}{{\small $z$}}
\psfrag{c}{{\small $z'$}}
\begin{center}
    \includegraphics[height=5cm]{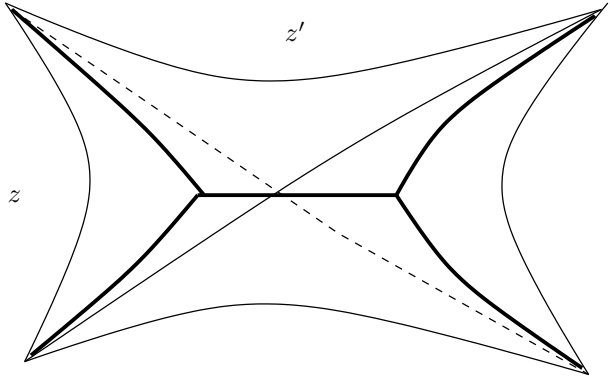}
\end{center}
  \caption{The geodesic spine}
  \label{fig:geodesic_spine}
\end{figure}

Let $\sigma$ be an ideal (topological) tetrahedron in $\widetilde{M},$ with shape parameters $z,z',z'',$ such that $z(Z_i)\to 1$ as $Z_{i} \to \xi.$  Regard $\sigma(Z_i)$ as an abstract metric space; first with the hyperbolic metric, then with this metric suitably rescaled.  The labels $z,z',z''$ are used for the edges of $\sigma(Z_i)$ without reference to their specific values at $Z_{i}.$ For each edge $e$ with parameter $z,$ join the centres of mass of the two faces containing $e$ by a geodesic arc and consider the common perpendicular to the $z$--edges. The rotational symmetry about this perpendicular interchanges the centres of mass of the two faces meeting in $e,$ and hence intersects the geodesic between them. Let the portion of the perpendicular between the two intersection points be the axis of the geometric spine, and add geodesic half--lines going from the endpoints of the axis to the vertices as indicated in Figure \ref{fig:geodesic_spine}. Note that this spine is well--defined for any value of $z,$ and denote it by $S_i(\sigma).$

Now assume that $\sigma$ is an ideal simplex in $\widetilde{M}$ which does not degenerate. Then let $S_i(\sigma)$ consist of the four geodesic arcs in $\sigma(Z_i)$ going from its centroid to its ideal vertices.


\subsection{Gromov--Hausdorff limit}

The following notions can be found in \cite{bh}, Chapter I.5. A subset $A$ of
a metric space $X$ is said to be $\varepsilon$--dense if every point of $X$ lies
in the $\varepsilon$--neighbourhood of $A.$ An $\varepsilon$--relation between two
metric spaces $X_1$ and $X_2$ is a subset $R \subseteq X_{1}\times X_{2}$ such
that the projection of $R$ onto each factor is $\varepsilon$--dense, and such
that if $(x_{1},x_{2}), (y_{1},y_{2}) \in R,$ then
\begin{equation}
| d_{X_{1}}(x_{1},y_{1}) - d_{X_{2}}(x_{2},y_{2}) | \le \varepsilon.
\end{equation}
The notation $X_{1} \sim_{\varepsilon} X_{2}$ indicates that there is an $\varepsilon$--relation
between them. If $X_{1} \sim_{\varepsilon} X_{2},$ then there is a
$3\varepsilon$--relation whose projection to each factor is onto. A sequence $\{
X_{i} \}$ of metric spaces converges to a metric space $X$ in the
Gromov--Hausdorff sense if there is a sequence $\{ \varepsilon_i \}$ of
non--negative real numbers such that $X_{i} \sim_{\varepsilon_{i}} X,$ and
$\lim_{i\to \infty} \varepsilon_i = 0.$


\subsection{The limiting tree}
\label{GH-limit of tet}

Let $v$ be an ideal vertex of an ideal triangle in $\H^3.$ There is a (geodesic) half--line perpendicular to the edge opposite $v$ which terminates in $v.$ The three perpendiculars thus obtained meet in a point inside the triangle which is called its \emph{centre of mass}. The \emph{geometric spine} of an ideal triangle is defined to be the union of the rays on the perpendiculars going from the centre of mass to the ideal vertices. Using the fact that all ideal triangles are congruent and that every hyperbolic triangle is contained in an ideal triangle, one can show:

\begin{lem}
The distance between the centre of mass in an ideal triangle in $\H^3$ and any of its edges is equal to $\ln \sqrt{3}$; the geometric spine of an ideal triangle is therefore $(\ln \sqrt{3})$--dense. In particular, every point on a side of a triangle is within distance $\ln 3$ of at least one point on the other two sides.
\end{lem}

\begin{figure}[t]
\begin{center}
    \includegraphics[width=9cm]{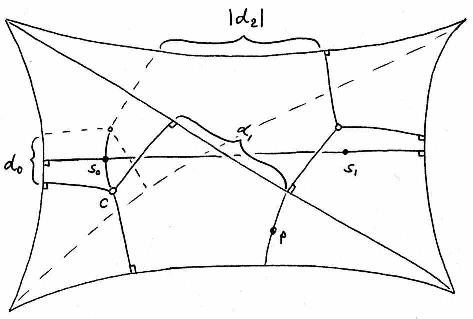}
\end{center}
  \caption{The geodesic spine is $\ln 27$--dense if $d_0 \le \ln 3$}
  \label{fig:spine dense}
\end{figure}

\begin{lem}\label{lem:spine dense}
Let $\sigma$ be an ideal (topological) tetrahedron in $\widetilde{M},$ with shape parameters $z,z',z'',$ such that $z(Z_i)\to 1$ as $Z_{i} \to \xi.$ If $|\ln |z(Z_i)|| \le \ln 3,$ then $S_i(\sigma)$ is $\ln 27$--dense in $\sigma(Z_i).$
\end{lem}

\begin{proof}
Since $\sigma(Z_i)$ is the convex hull of its vertices, it suffices to show that its faces are within distance $\ln 27$ from $S_i(\sigma).$ A face is in the $\ln \sqrt{3}$--neighbourhood of its spine, hence if the spines of all faces are in a $5 \ln \sqrt{3}$--neighbourhood of $S_i(\sigma),$ then the whole tetrahedron is in a $\ln 27$--neighbourhood of $S_i(\sigma).$ To simplify notation, $\sigma$ will be written instead of $\sigma(Z_i).$
  
Consider two faces of $\sigma$ meeting in an edge $e.$ On $e$ are the endpoints of the perpendiculars on the faces going into the vertices not in $e.$ The hyperbolic distance between them is precisely $|\ln |z(e)||.$ Since $zz'z''=-1,$ one has $\ln |z| + \ln |z'| + \ln |z''| = 0.$ Put $d_{j} = \ln |z^{(j)}|.$  For $z$ sufficiently close to one, $d_{1}> 0$ and $d_{2}< 0$ since $z' \to \infty$ and $z'' \to 0.$ Assume without loss of generality that $d_{0}\ge 0.$ Then $0 \le d_{0} < d_{1} \le | d_{2}| = d_{0} + d_{1},$ and by assumption $d_{0} \le \ln 3.$

Label the vertices of $\sigma$ by $v_{0},...,v_{3}$ such that $v_{0}$ and $v_{1}$ are the endpoints of an edge $e$ with parameter going to one. Consider the face opposite $v_3,$ and denote its centre of mass by $c.$ It needs to be shown that the geodesic rays $[c,v_i],$ $i=0,1,2,$ are within distance $5 \ln \sqrt{3}$ from $S_i(\sigma).$ Denote the singularities of $S_i(\sigma)$ by $s_{0}$ and $s_{1},$ such that $s_{0}$ is on the geodesic segment joining $c$ to its "opposite" centre of mass. Let $p$ be a point on $[c, v_2]$ which is at most distance $\ln \sqrt{3}$ from the endpoint of the perpendicular bisector through $v_3$ on a face containing the edge $[v_2, v_3].$ This point exists by  the previous lemma. The following inequalities can be read off from Figure \ref{fig:spine dense}:
\begin{align}
  d(c, s_0) &\le 2 \ln \sqrt{3} + d_{0}/2 \le 3 \ln \sqrt{3}\\
  d(p, s_{1}) &\le 4 \ln \sqrt{3} + d_{0}/2 \le 5 \ln \sqrt{3}\\
\label{eq:spines}  d(s_{0}, s_{1}) &\le \ln 27 + d_{0}+ d_{1}
\end{align}
This completes the proof for the rays $[c,v_i],$ $i=0,1$ by considering the hyperbolic triangles containing them and the point $s_0,$ and for $[p,v_2]$ by considering $[p,s_1,v_2].$ The endpoints of the geodesic segment $[c,p]$ are within distance $5 \ln \sqrt{3}$ from the geodesic segment $[s_0,s_1],$ and hence this must be true for the whole segment. This completes the proof for the particular face, and, by symmetry, for all faces.
\end{proof}

\begin{lem}\label{lem:global constant}
There is a constant $C(\xi)$ and $i_{1} \in \NN$ such that for each ideal tetrahedron $\sigma \subset \widetilde{M},$ $S_i(\sigma)$ is $C(\xi)$--dense in $\sigma(Z_i)$ for all $i \ge i_{1}.$ 
\end{lem}

\begin{proof}
As $Z_{i} \to \xi,$ we have $z_{j}(Z_{i}) \to 1$ for all $j = 1,...,k$ , and hence $\ln |z_{j}(Z_{i})| \to 0.$ Thus, there exists $i_{0},$ such that $|\ln z_{j}(Z_{i})| \le \ln 3$ for all $j = 1,...,k$ and for all $i \ge i_{0}.$
Lemma \ref{lem:spine dense} states that in this case $S_i(\sigma)$ is $\ln 27$--dense in each $\sigma$ which degenerates.

If all tetrahedra degenerate, then $i_{1}=i_{0}$ and $C(\xi) = \ln 27$ satisfy the requirements. Hence assume that $k<n.$  For each $j \in \{k+1,...,n\}$ choose an ideal hyperbolic tetrahedron of the limiting shape $\lim_{i\to\infty} z_{j}(Z_{i}),$ and construct the dual spine. Let $D$ be the maximum over the distances of the centroids to the respective singularities on the faces. Then for each limiting tetrahedron, the dual spine is $(D + \ln 3)$--dense. Let $E = D + \ln 9.$ Since the distances of the centroids to the singularities of faces vary continuously, there exists $i_{2}$ such that for all $z_{j},$ $j \in \{k+1,...,n\},$ and all $i \ge i_{2},$ the dual spine is $E$--dense.  Putting $C(\xi) = \max\{\ln 27, E\}$ and $i_{1} = \max\{i_{0},i_{2} \}$ completes the proof.
\end{proof}

Recall the definition of the dual spine. If the parameter of $\sigma$ is determined by the coordinate $z_j,$ then let $S(\sigma)$ be an abstract $\RR$--tree consisting of a closed interval of length $\xi_{3j-1}$ with two copies of $[0,\infty)$ attached at each of its endpoints. 

\begin{figure}[t]
\psfrag{b}{{\small $z$}}
\psfrag{c}{{\small $z'$}}
\psfrag{z}{{\small $z$}}
\psfrag{C}{{\small $(\sigma(Z_i), d)$}}
\psfrag{B}{{\small $(\sigma(Z_i), u(Z_{i})d)$}}
\psfrag{A}{{\small $S(\sigma)$}}
\psfrag{a}{{\small $l$}}
\begin{center}
    \includegraphics[width=14cm]{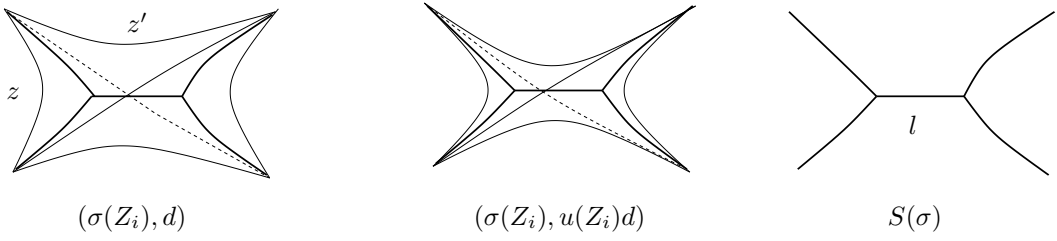}
\end{center}
  \caption{The Gromov-Hausdorff limit; $l=\lim_{i\to\infty} u(Z_i) \ln | z'(Z_i)|$}
  \label{fig:geom_spine}
\end{figure}

\begin{lem}\label{lem:spine is GH limit of tet}
There is a null--sequence $\{ \varepsilon_i\}$ of positive real numbers such that for each $\varepsilon_i$ there is a $\varepsilon_i$--relation $R \subset (\sigma(Z_i), u(Z_{i})d) \times S(\sigma)$ for each ideal 3--simplex $\sigma \subset \widetilde{M}.$ 
Thus, $S(\sigma)$ is a Gromov-Hausdorff limit of $\{(\sigma(Z_i), u(Z_i)d)\}.$
\end{lem}

\begin{proof}
Let $j \in\{1,...,k\},$ and $\sigma \subset \widetilde{M}$ a simplex with corresponding parameter triple $(z_{j}, z'_{j}, z''_{j}).$ Consider the sequence $(\sigma(Z_i), u(Z_{i})d),$ where $z_{j}(Z_{i}) \to 1.$  Equation (\ref{eq:spines}) can be written as 
$$d(s_{0}, s_{1}) \le \ln 27 + |\ln |z_{j}|| + \ln |z'_{j}|.$$ 
Hence $\lim_{i\to\infty} u(Z_{i})d(s_{0}(Z_{i}), s_{1}(Z_{i})) = \xi_{3j-1},$ which by assumption is a positive real number. After possibly passing to a subsequence of $\{Z_i\},$ one may therefore choose a null--sequence $\{\varepsilon_i\}$ such that $|\xi_{3j-1} - u(Z_i)d(s_{0}(Z_{i}), s_{1}(Z_{i}))| \le \varepsilon_{i}.$  
Letting $\delta_i = \max\{\varepsilon_i, u(Z_{i})C(\xi)\},$ one has $\lim_{i\to\infty}\delta_i = 0.$

Label the singularities of $S(\sigma)$ by $S_{0}$ and $S_{1},$ and its four ideal vertices by $V_{0},$...,$V_{1}.$ Denote the distance function on $S(\sigma)$ by $d_{S}.$ A relation $R \subset \sigma(Z_i) \times S(\sigma)$ is defined as follows. Label the ideal vertices and the singularities of the dual spine of $\sigma(Z_i)$ with reference to $\sigma,$ so that they are chosen consistently for all $i.$ Then define $R$ between $G_i(\sigma)$ and $S(\sigma)$ by including $(s_i,S_i),$ and extending the relation isometrically on the infinite ends, and by scaling linearly on the intervals between the singularities.

The projection of $R$ to the second factor is onto, and the projection to the first factor is $u_{i}C(\xi)$--dense, and in particular $\delta_{i}$--dense. Since the ends are identified isometrically, the maximal length distortion occurs on paths containing both singularities. One has:
\begin{align*}
  & |u(Z_{i})d(s_{0}(Z_{i}), s_{1}(Z_{i})) - d_{S}(S_{0}, S_{1})|\\
  =  &  |u(Z_{i})d(s_{0}(Z_{i}), s_{1}(Z_{i})) - \xi_{3j-1}| \\
  \le & \delta_i.
\end{align*}
Whence $R$ is a $\delta_{i}$--relation between $(\sigma(Z_i), u(Z_{i})d)$ and $S(\sigma).$ A similar argument applies to the tetrahedra which do not degenerate.
\end{proof}

\begin{pro}
The sequence $(\widetilde{M}_{i}, u(Z_{i})d)$ converges to the leaf--space $\widetilde{M}/\widetilde{\fol}$ of the transversely measured singular codimension--one foliation in the sense that:
\begin{equation*}
\lim_{i\to \infty} \{ (\sigma(Z_i), u(Z_i)d) : \sigma \subset \widetilde{M} \} / \Psi_i \quad = \quad \{ S(\sigma) : \sigma \subset \widetilde{M} \} / \sim.
\end{equation*}
\end{pro}

\begin{proof}
To show that the proposition follows from Lemma \ref{lem:spine is GH limit of tet}, one needs to show that the limiting "face pairings" amongst the dual spines give the equivalence relation described earlier. This follows from the following two observations. The limit of each face is of the shape of a $Y,$ and since the gluings of faces are by isometries, two $Y$ shapes are identified isometrically along rays corresponding to the edges of tetrahedra. An isometric gluing identifies the centres of mass of the faces, which correspond to the vertices of the $Y$ shapes. Hence the singularities glue up, and the identification is isometrically along the three infinite rays.
\end{proof}

The above does not show that the sequence $\{(\widetilde{M}_{i}, u(Z_{i})d)\}$ converges in the Gromov--Hausdorff sense to $\widetilde{M}/\widetilde{\fol}$ or $\tree_N$ (though this is true for arbitrarily large simplicial subsets), but it gives a useful device for comparing the action on $\tree_N$ with the limiting action on $\H^3.$ This is similar to the situation in \cite{best}.


\subsection{Translation length functions}
\label{finish}

This section continues with the notation introduced in Section \ref{sec:Outline and definitions}. For each $\gamma \in \pi_1(M),$ consider the following two translation lengths:
\begin{align*}
l_N(\gamma) &= \inf \{ d_N(x,\gamma x) : x \in \tree_N\}, \\
l_\H(\rho_i(\gamma)) &= \inf \{ d(x,\rho_i(\gamma) x) : x \in \H^3 \},
\end{align*}
where $d_N$ denotes the distance function on $\tree_N.$ 

\begin{pro}\label{pro: length functions}
For each $\gamma \in \pi_1(M),$ the sequence $\{ u(Z_i) l_\H(\rho_i(\gamma)) \}$ contains a convergent subsequence, and its limit $l_\xi(\gamma)$ satisfies:
\begin{equation}\label{eq: length functions ineq}
0 \le l_\xi(\gamma) \le l_N(\gamma).
\end{equation}
\end{pro}

\begin{proof}
Let $x \in \widetilde{M},$ and $\tilde{\gamma}$ be a path from $x$ to $\gamma x$ in $\widetilde{M}.$ Then:
\begin{equation*}
0 \le l_\H(\rho_i(\gamma))  \le d_\H (D_i(x), \rho_i(\gamma)D_i(x)) = 
d_\H (D_i(x), D_i(\gamma x)) \le  \text{length}(D_i(\tilde{\gamma})).
\end{equation*}
The map $\widetilde{M}_i \to D_i(\widetilde{M})$ restricts to an isometry on each ideal tetrahedron, so
$\text{length}(D_i(\tilde{\gamma})) =  \text{length}(f_i(\tilde{\gamma})),$ and one has
\begin{equation}\label{eq: proof of transl length}
0 \le u(Z_i) l_\H(\rho_i(\gamma))  \le u(Z_i) \text{length}(f_i(\tilde{\gamma})).
\end{equation}

Recall from \cite{part1} that a path is called admissible if it consists of finitely many sub--paths each of which is either contained in a leaf or is transverse to the foliation, and that the distance between two points in $\tree_N$ is the infimum over the measures of all paths connecting them; the transverse measure is denoted by $\mu.$ Any path can be deformed into an admissible path without increasing the measure.

Let $p \in \tree_n,$ and assume that $x \in \widetilde{M}$ is a preimage of $p.$ Then for any $\varepsilon>0$ there is an admissible path $\tilde{\gamma}$ in $\widetilde{M}$ joining $x$ and $\gamma x,$ such that
\begin{equation*}
d_N(p,\gamma p) \le \mu (\tilde{\gamma}) + \varepsilon.
\end{equation*}
The path can be subdivided into $t$ sub--paths for some $t \in \NN,$ such that each of them is contained in an ideal tetrahedron.

According to Lemma \ref{lem:spine is GH limit of tet}, for some $n\in \NN$ and all $i>n$ and all $\sigma \subset \widetilde{M},$ there is a $\varepsilon_i$--relation between $\sigma(Z_i)$ and $S(\sigma)$ with $\varepsilon_i< \varepsilon.$ This relation extends to a surjective $(3\varepsilon_i)$--relation, and the properties of the Gromov--Hausdorff limit imply that:
\begin{equation*}
 \mid \mu(\tilde{\gamma})- u(Z_i)\text{length}(f_i(\tilde{\gamma}))\mid \le 3t \varepsilon,
\end{equation*}
giving
\begin{equation*}
 \mid d_N(p,\gamma p)  - u(Z_i)\text{length}(f_i(\tilde{\gamma}))\mid \le (3t+1) \varepsilon.
\end{equation*}
Since $\varepsilon>0$ is arbitrary, inequality (\ref{eq: proof of transl length}) implies:
\begin{equation*}
d_N(p,\gamma p)  \ge u(Z_i) \inf_{\tilde{\gamma}} \{ \text{length}(f_i(\tilde{\gamma})) \} \ge u(Z_i) l_\H(\rho_i(\gamma)) \ge 0,
\end{equation*}
where the infimum is taken over all paths from $x$ to $\gamma x.$ Hence the sequence $\{ u(Z_i) l_\H(\rho_i(\gamma)) \}$ is bounded, and since $p \in \tree_N$ is arbitrary, the conclusion follows.
\end{proof}

The above proposition gives a sufficient condition ($l_\xi \neq 0$) for the action on $\tree_N$ to be non-trivial, and a necessary condition ($l_\xi = 0$) for the action on $\tree_N$ to be trivial. This is first translated into a condition only involving information from spun-normal surface theory, and then into a general statement linking the ideal points of the deformation variety to the compactification of the character variety due to Morgan and Shalen.

\begin{cor}\label{pro: non-trivial action via peripheral}
Let $\xi \in \D_\infty(\tri)$ and $N=N(\xi).$ Assume that for some peripheral element $\gamma \in \pi_1(M),$ $\nu_N(\gamma)\neq 0.$ Then the action of $\pi_1(M)$ on $\tree_N$ is non-trivial.
\end{cor}

\begin{proof}
Let $\gamma \in \pi_1(M).$ The hyperbolic translation length $\l_i$ of $\rho_i (\gamma)$ satisfies $\tr \rho_i(\gamma) = e^{\l_i/2} + e^{-\l_i/2}.$ The relationship between the growth rates of parameters and the valuation of squares of eigenvalues of peripheral elements described in the proof of Proposition \ref{comb:essential prop} gives that $\nu_N(\gamma)\neq 0$ implies $l_\xi(\gamma) > 0.$ The result now follows from Proposition \ref{pro: length functions}.
\end{proof}

\begin{cor}\label{cor: non-trivial action 2}
Let $\xi \in \D_\infty(\tri),$ $N=N(\xi),$ and $(Z_i )\subset \D(\tri)$ be a sequence strongly converging to $\xi.$ The action of $\pi_1(M)$ on $\tree_N$ is non-trivial, if an ideal point of the character variety is approached by the sequence $\chi_\tri(Z_i)$
\end{cor}

\begin{proof}
If the sequence $\chi_\tri(Z_i)$ approaches an ideal point of the character variety, i.e.\thinspace some trace becomes unbounded, we need to show that the scaling factors $u(Z_i)$ are not too big as to render $l_\xi$ trivial. This follows from the facts that $\chi_\tri(Z_i)$ is a rational function in the shape parameters, and that rescaling any shape parameter tending to $0$ or $\infty$ by $u(Z_i)$ gives a non-zero number proportional to its growth rate.
\end{proof}

\begin{rem}
It follows from the construction that if the length function $l_\xi$ is non-trivial, it determines the point in the Morgan-Shalen compactification with coordinate $(l_\xi(\gamma))_{[\gamma] \in \mathcal{C}},$ where $\mathcal{C}$ is the set of conjugacy classes of elements in $\pi_1(M).$ 
\end{rem}


\subsection{Dual surfaces}
\label{sec:limiting surface}

Throughout this section, assume that $\xi$ has rational coordinate ratios, and recall that $S(\xi)$ is a 2--sided surface. In this case $\tree_N = \widetilde{M}/\widetilde{\fol},$ and it can be rescaled so that it contains the \emph{dual (simplicial) tree} $\tree_S$ of the lift of $S=S(\xi)$ to $\widetilde{M}$ as a subtree. A surface is \emph{non-trivial} if it is essential or can be reduced to an essential surface by performing compressions and then possibly discarding some components. Each vertex stabiliser of the action on $\tree_S$ is conjugate to the image under the inclusion map $\pi_1(M_i)\to\pi_1(M)$ for some component $M_i$ of $M-S,$ and each edge stabiliser to the image of $\pi_1(S_j)\to\pi_1(M)$ for some component $S_j$ of $S.$ Thus, $S$ is non-trivial if and only if none of the images $\pi_1(M_i)\to\pi_1(M)$ is onto, and this is the case if and only if the action of $\pi_{1}(M)$ on $\tree_N$  is non-trivial. 

\begin{pro}\label{pro:dual surfaces}
Let $\{Z_i\}$ be a sequence strongly converging to $\xi,$ and $\chi_i = \varphi_\tri(Z_i).$
\begin{enumerate}
\item The sequence $\{ \chi_i|C\} $ converges in $\PX(C)$ for each connected component $C$ of $M-S.$
\item The sequence $\{ \chi_i|{S_j}\}$ converges in $\PX(S_j)$ to a reducible character for each connected component $S_j$ of $S.$
\end{enumerate}
In particular, if for some $\gamma\in \pi_1(M),$ $\{ | \chi_i(\gamma)| \}$ is unbounded, then $S$ is non-trivial and (weakly) dual to an ideal point of a curve in the character variety of $M$; hence there is a $\pi_1(M)$--equivariant map from $\tree_N$ to the Bass--Serre tree associated to the ideal point.
\end{pro}

\begin{proof}
Corollary 2 in \cite{serre}, I6.2, states that a finitely generated group acting on a simplicial tree fixes a vertex if and only if all generators and their double products have zero translation length. This together with equation (\ref{eq: length functions ineq}) implies that $\{ \chi_i|C\}$ is contained in a compact subset of $\PX(C).$ To prove that there is a unique accumulation point, some terminology and observations are needed.

The spun-normal surface $S$ divides each tetrahedron into several types of regions. A region which contains two normal isotopic discs is a trivial $I$--bundle over a normal disc and will be called a \emph{slab}. A region which is not a slab is a \emph{vertex region} if it is a corner cut off by a triangle, otherwise it is a \emph{thick region}. A thick region is a truncated tetrahedron or a truncated prism.

Assume that $C$ contains some thick region. This region meets some face of an ideal tetrahedron $\sigma$ in $M$ in a truncated ideal triangle. Choose a lift of this face to $\widetilde{M}$ and denote it by $F.$ Define the maps $D_i$ such that $F$ is sent to the ideal triangle $[0,1,\infty]$ for all $i.$ Each region is associated to a vertex on a dual spine as in Section \ref{subsec:dual tree}. Let $p$ be the vertex corresponding to the centre of mass of $F.$ Let $\gamma$ be an element of the group $\Gamma$ in the conjugacy class of the image $\im(\pi_1(C) \to \pi_1(M))$ which stabilises $\tilde{C}\subset \widetilde{M}.$ Then $\gamma(p)$ is contained in $\tilde{C}$ and hence within bounded hyperbolic distance of the vertex corresponding to another region of $\tilde{C},$ say $p'.$ There is a finite sequence of adjacent regions and their vertices $p=p_0,{\ldots} ,p_k=p'$ such that the regions corresponding to $p_j$ and $p_{j+1}$ are identified along a face of some tetrahedron in $\widetilde{M}.$  As the tetrahedra degenerate, the images of $p_j$ and $p_{j+1}$ under $D_i$ stay within bounded hyperbolic distance for all $i,$ and hence $p$ and $p'$ stay within bounded hyperbolic distance. Thus, $\lim_{i\to \infty} u(Z_i)d(D_i(p),D_i(p'))=0,$ and the sequence $\{ \rho_i(\gamma)\}$ converges since it only depends on the limiting shapes of the tetrahedra and the combinatorial gluing data. Since $\gamma$ was chosen arbitrarily, the sequence $\{ \rho_i|_C\}$ converges algebraically, and the conclusion follows since choosing a different developing map changes the representation by conjugation. This proves (1) for components that contain thick regions.

A normal disc in $S$ is called \emph{bad} if it is contained in the boundary of some thick region. If $S$ contains no normally isotopic components, then each connected component of $S$ meets some thick region in a bad disc. Thus, if $S_j$ does not contain a bad disc, one may replace it by a normally isotopic component for the purpose of proving (2). Associated with the thick region met by $S_j$ is a centre of mass of some face in $\tri.$ Define the maps $D_i$ such that a fixed lift of this face is sent to the ideal triangle $[0,1,\infty]$ for all $i.$ Similar to the above, using the edges in dual spines associated to $S_j,$ one sees that $\{ \rho_i|_{S_j}\}$ converges algebraically. Since the limiting representation of $\pi_1(S_j)$ fixes a vertex of $[0,1,\infty],$ it is reducible. Whence (2) holds.

To conclude the proof of (1), suppose that the component $C$ contains no thick region. If it contains a vertex region, then it must be entirely made up of vertex regions, so its boundary is a vertex linking surface; $S$ contains no such component.
If it contains a slab region, then it must be entirely made up of slabs, in which case it is either a trivial or a twisted $I$--bundle over a spun-normal surface, and $S_j \subseteq \partial C$ for some $j.$ Both cases follow from (2); the first trivially, the second from the fact that $\chi_i|C$ is uniquely determined by $\chi_i|{S_j}$ since $[\pi_1(C):\pi_1(S_j)]=2$ and $(\tr A)^2 = \tr (A^2)-2$ for each $A \in \SL.$

The last part of the proposition follows directly from \cite{tillus_ei}, Lemma~6, and \cite{tillus_mut}, Section~3. The former shows that the sequence $\{Z_i\}$ may be chosen on a curve in $\D(\tri).$ The sufficient condition ($l_\xi \neq 0$) for the action to be non-trivial, together with (1) and (2), allows the construction in \cite{tillus_mut} to be applied to yield the assertion. Note that the map to the Bass--Serre tree may not be an isometry (even up to scaling and restriction to maximal invariant subtrees).
\end{proof}

\begin{rem}
If the above necessary condition for a trivial action ($l_\xi=0$) is not sufficient, then there may be closed essential surfaces which are detected by the deformation variety, but not by the character variety.
\end{rem}

The proofs of the main results are completed with the following immediate consequence of Corollary \ref{pro: non-trivial action via peripheral} and Proposition \ref{pro:dual surfaces}.

\begin{cor}\label{cor:thm 1 part 2}
If $S(\xi)$ is spun--normal, then $S(\xi)$ is non-trivial and (weakly) dual to an ideal point of a curve in the character variety.
\end{cor}


\section{The figure eight knot}
\label{sec: figure eight}

\begin{figure}[t]
\psfrag{a}{{\small $z$}}
\psfrag{b}{{\small $z'$}}
\psfrag{c}{{\small $z''$}}
\psfrag{w}{{\small $w$}}
\psfrag{x}{{\small $w'$}}
\psfrag{y}{{\small $w''$}}
\psfrag{z}{{\small $z$}}
\psfrag{0}{{\small $0$}}
\psfrag{1}{{\small $1$}}
\psfrag{i}{{\small $\infty$}}
\psfrag{h}{{\small $zw$}}
    \begin{center}
      \includegraphics[width=9cm]{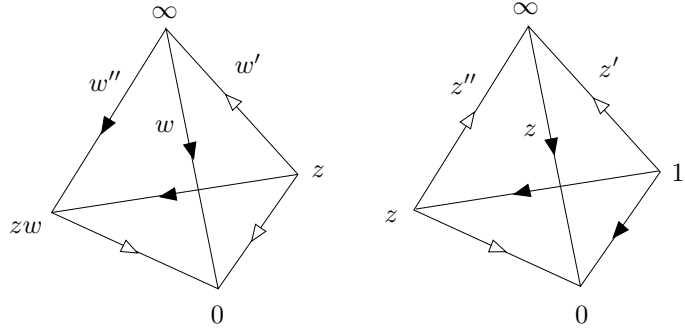}
    \end{center}
    \caption{An ideal triangulation of the figure eight knot complement}
    \label{fig:fig8_tets}
\end{figure}

Let $M$ denote the complement of the figure eight knot. An ideal triangulation of $M$ is shown in Figure \ref{fig:fig8_tets}. The quadrilateral types dual to $w^{(k)}$ and $z^{(k)}$ will be denoted by $p^{(k)}$ and $q^{(k)}$ respectively.

One has the following hyperbolic gluing equations:
\begin{align*}
1 &= (w')^2 w'' (z')^2 z'', \\
1 &= w^2 w'' z^2 z''.
\end{align*}
This gives the matrix
\begin{equation*}
A = \begin{pmatrix} 0 & 2 & 1 & 0 & 2 & 1 \\
                      2 & 0 & 1 & 2 & 0 & 1 \end{pmatrix}.
\end{equation*}
Hence
\begin{equation*}
B=A C_2 = \begin{pmatrix} -1 & -1 & 2 & -1 & -1 & 2 \\
                          1 & 1 & -2 & 1 & 1 & -2 \end{pmatrix},
\end{equation*}
which determines a single $Q$--matching equation:
\begin{equation*}
0 = q + q' - 2 q'' + p + p' - 2 p''.
\end{equation*}
This agrees with \cite{part1}. One can also work out the induced triangulation of the torus end and determine standard generators for the peripheral subgroup, giving $\nu (\l ) =  - 2 q - 2 q' + 4 q''$ and
$\nu (\m ) = -p'+p''-q+q''.$ The set $\N (\tri)$ is $0$--dimensional, its elements are scaled to integer solutions and listed in Table \ref{tab:surfaces fig8}. All spun-normal surfaces corresponding to these solutions are once--punctured Klein bottles. Since $\D (\tri) \neq \emptyset,$ its follows from the symmetries that $\D_\infty(\tri)= \N (\tri).$

\begin{table}[h]
\begin{center}
\begin{tabular}{ c | r | r | r}
solution & $\nu (\m )$  & $\nu (\l )$ & slope \\
\hline
(2,0,0,0,0,1) & 1 &4  &        --4  \\
(0,2,0,0,0,1) & --1 & 4  &       4   \\
(0,0,1,2,0,0) & --1 & --4  & --4 \\
(0,0,1,0,2,0) &  1&  --4  &     4 \\
\end{tabular}
\end{center}
\caption{Normal surfaces in the figure eight knot complement}
\label{tab:surfaces fig8}
\end{table}


\subsection{Face pairings}

A presentation of the fundamental group can be worked out from the triangulation, and simplifies to:
\begin{equation}\label{fig8:fund2}
    \pi_1(M) = \langle A, \m \mid A\m A^{-2}\m = \m A\m A^{-1}\rangle.
\end{equation}
Furthermore, $\m$ is a meridian, and the corresponding longitude is
\begin{equation*}
\l = \m^{-1}A\m A^{-1} \m^{-1}A^{-1}\m A.
\end{equation*}
Let a fundamental domain for $M$ be given by embedding the tetrahedra in $\H^3$ as indicated by the vertex labels in Figure \ref{fig:fig8_tets}. The face pairings associated to $Z = (w,w',w'',z,z' ,z'') \in \D (\tri)$
are defined by assignments of the following ordered triples:
\begin{align*}
A_Z &: [0,z,zw] \to [\infty, 1,z]\\
\m_Z &: [\infty, 0,zw] \to [1,0,z]\\
G_Z &: [\infty, zw,z] \to [\infty, 0,1]
\end{align*}
It follows that $G_Z = A_Z\m_Z A_Z^{-1}.$  From above face pairings and the equation $z(z-1)w(w-1)=1,$ one obtains representations into $\PSL$ by putting
\begin{align*}
\rho_Z(\m )&= \frac{1}{\sqrt{w(1-z)}} \begin{pmatrix} 1 & 0\\
                      1 & w(1-z) \end{pmatrix}\\
\rho_Z (A) &= \begin{pmatrix} 1 - wz & -(1-w)^{-1}\\
                      1-w & 0 \end{pmatrix}\\
\rho_Z(\l )&= \begin{pmatrix} w(w-1) & 0\\
                      z(1+w-w^2)(wz-w-z) & z(z-1) \end{pmatrix}
\end{align*}
This is in fact a representaton into $\SL,$ and the lower right entries in $\m$ and $\l$ correspond to square roots of the holonomies given by Thurston in \cite{t}, where $\mu(\m ) = w(1-z)$ and $\mu(\l )= z^2(z-1)^2.$ The image of $\D (\tri)$ in the $\PSL$--eigenvalue variety is parameterised by:
\begin{equation*}
m^4-2m^3-3m^2+2m-l+6-l^{-1}+2m^{-1}-3m^{-2}-2m^{-3}+m^{-4}=0,
\end{equation*}
where $m=\mu(\m )$ and $l=\mu(\l ).$ The map from $\D (\tri)$ to the $\PSL$--character variety can also be determined from the above. One has $(\tr \m)^2 = w + 2(1-wz)+z$ and $\tr A = 1 - wz.$ Putting $(\tr \m)^2=X$ and $\tr A = y,$ the image of $\D (\tri)$ in the $\PSL$--character variety is parameterised by:
\begin{equation*}
1 - y - y^2 + (y-1)X = 0,
\end{equation*}
giving a sphere in $\C P^2.$ The smooth projective model of $\D (\tri)$ is a torus, and it can be verified that the map $\D (\tri) \to \PX_0 (M)$ is generically 2--to--1.


\subsection{Limiting characters at ideal points}

Given the symmetries of $M$ and its the triangulation, it suffices to consider a degeneration of the ideal triangulation to one of the four ideal points of $\D (\tri).$ The point whose associated spun-normal surface coordinates are $(0,2,0,0,0,1)$ will be chosen, and ``geometric'' degenerations, i.e. degenerations where both tetrahedra stay positively oriented and only in the limit become degenerate will be studied. The deformation variety is birationally equivalent to the variety in $\C^2$ defined by the single equation:
\begin{equation*}
   z(1-z)w(1-w)=1.
\end{equation*}
Solving the above equation in terms of $w,$ gives
\begin{equation} \label{fig8:solve}
   z = \frac{1}{2}\Big(1 \pm \sqrt{1 + \frac{4}{w(w-1)}}\Big)
\end{equation}
At the complete structure, one has $w_0=z_0=\frac{1}{2}(1 + \sqrt{-3})$; hence take the solution for $z$ with positive sign in front of the root. The desired ideal point corresponds to the degeneration
$w \to 0.$ Note that then $z\to \infty$ at half the rate. Let $w = w(r) = r^2 w_0,$ and obtain a power series expansion of $z$ for $r$ around zero using equation (\ref{fig8:solve}). This is of the form $z(r) = \frac{w_0}{r} + \frac{1}{2} + \varphi(r),$ where $\varphi(0) = 0.$

All points on the path $[1,0) \to \D (\tri)$ given by $r \to (w(r),z(r))$ correspond to geometric solutions to the hyperbolic gluing equations (see \cite{t}). The face pairings can be used to determine the limiting representations. As $r \to 0,$ one has the following limiting traces:
\begin{equation*}
    \tr \rho_Z(A) \to 1 \qquad\text{and}\qquad \tr^2 \rho_Z(\m) \to \infty,
\end{equation*}
whilst the eigenvalue of $\m^4 \l,$ which is equal to $-w^2z(1-z)^3,$ approaches one. Thus, $\m^4 \l$ is a strongly detected boundary slope.

The limiting splitting of $M$ corresponds to a splitting of $M(4,1),$ which is a graph manifold, along an essential torus (see \cite{chk}). The limiting pieces admit Seifert fibered structures, and are a twisted $I$--bundle over the Klein bottle, $S_1,$ and a trefoil knot complement, $S_2.$ Note that in this case passing to an orientable surface introduces a certain redundancy; splitting along the Klein bottle gives a nicer decomposition.

\begin{figure}[t]
\psfrag{a}{{\small $A$}}
\psfrag{m}{{\small $\m$}}
\psfrag{1}{{\small $k_1$}}
\psfrag{2}{{\small $k_2$}}
    \begin{center}
     \includegraphics[width=7cm]{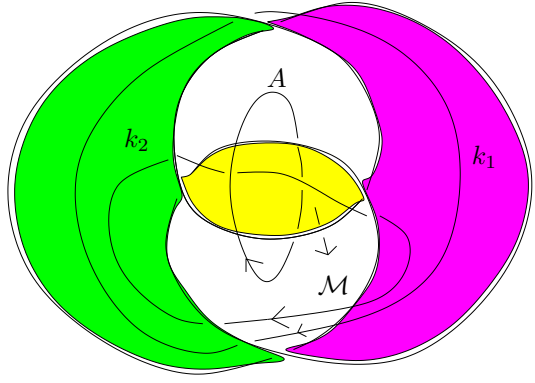}
    \end{center}
    \caption{The punctured Klein bottle}
    \label{fig:fig8_klein}
\end{figure}

The fundamental groups of the complementary pieces can be worked out from Figure \ref{fig:fig8_klein}.
Let $K'$ be the punctured Klein bottle in $M$ shown in Figure \ref{fig:fig8_klein}, and let $K$ be the corresponding Klein bottle in $M(4,1).$ Identify $I \tilde{\times} K$ with a regular neighbourhood of $K$ in $M(4,1).$ Standard generators for $\im (\pi_1(K) \to \pi_1(M(4,1)))$ are
$k_1 = \m A^{-2} \m$ and $k_2 = A^{-1} \m A^{-1}.$ One has
\begin{equation*}
k_2k_1k_2^{-1}k_1 = A \m^{-1} A^{-1} (\m^4 \l ) A\m A^{-1} = 1.
\end{equation*}
Standard generators for the boundary torus of $I \tilde{\times} K$ are $\m_1 = k_1$ and $\l_1 = k_2^2.$

Generators for the complement of $I \tilde{\times} K$ in $M(4,1)$ are given by $u = A$ and $v = A\m A\m A^{-1}.$ One has
\begin{equation*}
v u^{-3} v = A \m (\m^4 \l ) \m^{-1} A^{-1} =1.
\end{equation*}
Generators for the boundary torus of the trefoil knot complement are $\m_2 = u^{-1} v$ and $\l_2 = u^3$ (the meridian is standard, but the longitude is not).

The decomposition of the fundamental group of $M(4,1)$ can be worked out from the above. It is an amalgamated product of $\pi_1(S_1) = \langle k_1,k_2 \mid k_2k_1k_2^{-1}k_1 = 1\rangle$ and $\pi_1(S_2) = \langle u,v \mid u^3 = v^2\rangle,$ amalgamated by $\m_1 = \m_2^{-1}$ and $\l_1 = \l_2^{-1} \m_2.$

The limiting representation on the trefoil knot complement is determined by:
\begin{equation*}
\rho (A) \to \begin{pmatrix} 1 & -1 \\ 1 & 0 \end{pmatrix},
\rho (A \m A\m A^{-1}) \to \begin{pmatrix} 0 & -1 \\ 1 & 0 \end{pmatrix},
\rho (\m A\m A^{-1}) \to \begin{pmatrix} 1 & 0 \\ 1 & 1 \end{pmatrix}.
\end{equation*}
This representation corresponds to a 2--dimensional hyperbolic structure on the base orbifold of $S_2,$ which is a $(2,3,\infty )$--turnover. The limiting representation of the boundary torus is:
\begin{equation*}
\rho( \m_2 ) \to \begin{pmatrix} 1 & 0 \\ 1 & 1 \end{pmatrix},
\qquad
\rho (\l_2 ) \to \begin{pmatrix} -1 & 0 \\ 0 & -1 \end{pmatrix}.
\end{equation*}

In order to obtain a (finite) limiting representation of $S_1,$ one needs to conjugate the face pairings by a diagonal matrix with eigenvalues $(r w_0^{-1})^{\pm 1/4}.$ The limiting representation is then:
\begin{equation*}
\rho (k_1) \to \begin{pmatrix} 1 & 0 \\ 0 & 1 \end{pmatrix},
\qquad
\rho(k_2) \to \begin{pmatrix} 0 & 1 \\ -1 & 0 \end{pmatrix},
\end{equation*}
giving a cyclic group of order two in $\PSL.$ The limiting image of $\rho( \m_1 )$ is $E,$ and the limiting image of $\rho( \l_1 )$ is $-E.$


\subsection{Description of the singular foliation}

\begin{figure}[h!]
\begin{center}
   \subfigure[$H_1$]{
        \includegraphics[height=3.7cm]{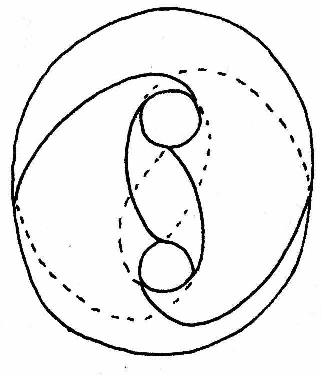}}\qquad
   \subfigure[$\fol \cap \partial B_1$]{
        \includegraphics[height=3.7cm]{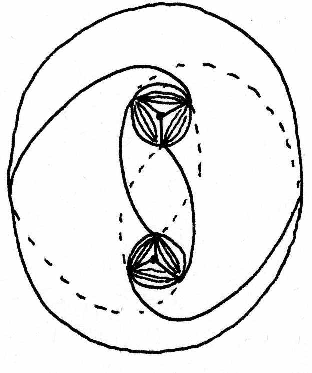}
        \label{fig:fol fig8 B_1}}
        \qquad
      \subfigure[]{
        \includegraphics[height=3.7cm]{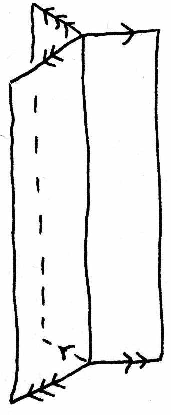}
        \label{fig:fol fig8 sing curve}}\qquad
      \subfigure[1--handle]{
        \includegraphics[height=3.7cm]{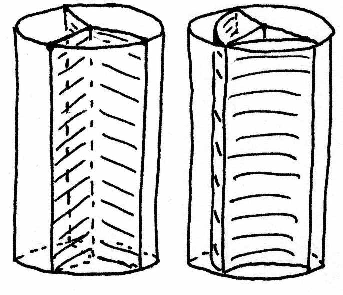}
        \label{fig:fol fig8 1-handle}}
        \\
      \subfigure[The singular leaf meeting $B_2$]{
        \includegraphics[height=4cm]{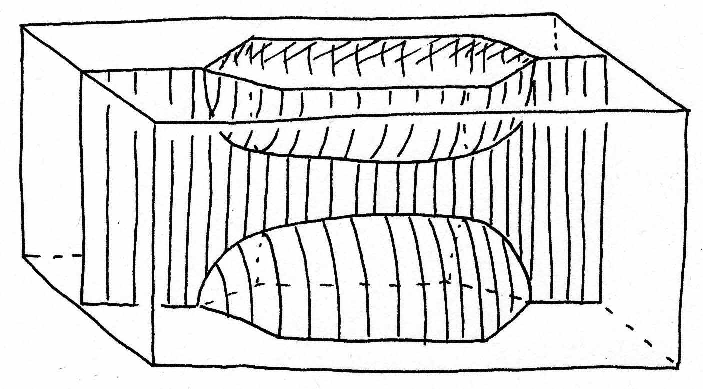}
         \label{fig:fol fig8 B_2 sing}}
      \subfigure[A non--singular leaf meeting $B_2,$ the thick parts of the boundary are attached to squares in the 1--handles, the thin parts run along the knot]{
        \includegraphics[height=4cm]{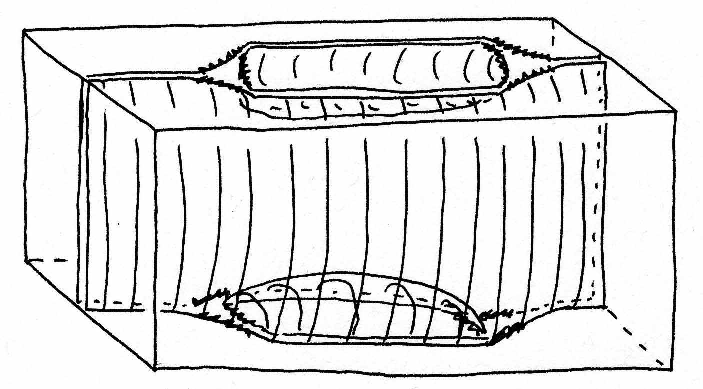}
         \label{fig:fol fig8 B_2 non--sing}}
   \end{center}
    \caption[]{Singular foliation of the figure eight knot complement at an ideal point}
    \label{fig:fol fig8}
\end{figure}

The topology of the leaves in the singular foliation $\fol$ associated to the above ideal point can be determined from the spun-normal surface data. One leaf is the once--punctured Klein bottle made up of three quadrilaterals and infinitely many triangles; the single singular leaf is made up of the four singular pieces plus one quadrilateral and infinitely many triangles -- removing the triple curve yields a once--punctured M\"obius band; all other leaves are twice punctured tori made up of six quadrilaterals and infinitely many triangles.

An explicit description of $\fol$ was found by Thurston \cite{t81}, and is given in the following commentary to Figure \ref{fig:fol fig8}. Thickening the surface shown in Figure \ref{fig:fig8_klein} yields a genus--two handlebody $H_1$ in $S^3$ with the figure eight knot lying on its boundary, and $H_1$ minus the knot is foliated by one once--punctured Klein bottle and parallel twice--punctured tori. It remains to describe the foliation of the complementary handlebody $H_2$ minus the knot.

A small neighbourhood of the singular curve on the singular leaf in $\fol$ is shown in Figure \ref{fig:fol fig8 sing curve}. The singular curve loops through the two holes of $H_1.$ Near these holes, the foliation is as shown in Figures \ref{fig:fol fig8 B_1} and \ref{fig:fol fig8 1-handle}. So attaching a small 1--handle to each hole of $H_1$ which meets $\fol$ in the shown product gives a ball $B_1$ in $S^3,$ with the figure eight knot contained partly on its boundary, partly in its interior. The foliation of the complementary ball $B_2$ minus the knot is shown in Figures \ref{fig:fol fig8 B_2 sing} and \ref{fig:fol fig8 B_2 non--sing}. There is one singular leaf consisting of two boats attached along their keels, and the remaining leaves are discs. Any non--singular leaf contained in $H_2$ meets $B_2$ in exactly one discs in each of the four components of $B_2$ minus the singular leaf, and $B_1$ in six rectangles; one in each component of the 1--handles minus the singular leaf. The foliation of $H_2$ minus the knot by one singular leaf and parallel twice--punctured tori is now obtained by appropriately attaching two foliated 1--handles to $B_2.$



\address{Stephan Tillmann\\ School of Mathematics and Physics\\ The University of Queensland\\ Brisbane, QLD 4072, Australia\\
(tillmann@maths.uq.edu.au)\\}
\Addresses


\end{document}